\documentclass[10pt]{amsart}

\usepackage{amsfonts, amstext, amsmath, amsthm, amscd, amssymb}
\usepackage{epsfig, graphics, psfrag}
\usepackage{overpic}
\usepackage{color}
\usepackage{our-macros}

\begin{document}

\title{Angled decompositions of arborescent link complements}
\date{\today}

\author{David Futer}
\thanks{Futer partially supported by NSF grant DMS-0353717 (RTG)}
\address{David Futer, Mathematics Department, Michigan State
  University, East Lansing, MI 48824, USA}
\email{dfuter@math.msu.edu}

\author{Fran\c{c}ois Gu\'eritaud}
\thanks{Gu\'eritaud partially supported by NSF grant DMS-0103511}
\address{Fran\c{c}ois Gu\'eritaud, DMA -- UMR 8553 (CNRS), \'Ecole normale 
sup\'erieure, 45 rue d'Ulm, 75005 Paris, France}
\email{Francois.Gueritaud@ens.fr}

\begin{abstract}
This paper describes a way to subdivide a 3--manifold into \emph{angled blocks}, namely polyhedral pieces that need not be simply connected. When the individual blocks carry dihedral angles that fit together in a consistent fashion, we prove that a manifold constructed from these blocks must be hyperbolic. The main application is a new proof of a classical, unpublished theorem of Bonahon and Siebenmann: that all arborescent links, except for three simple families of exceptions, have hyperbolic complements.
\end{abstract}

\maketitle

\section{Introduction}\label{sec:intro}

In the 1990s, Andrew Casson introduced a powerful technique for constructing and studying cusped hyperbolic 3--manifolds. His idea was to subdivide a manifold $M$ into \emph{angled ideal tetrahedra}: that is, tetrahedra whose vertices are removed and whose edges carry prescribed dihedral angles. When the dihedral angles of the tetrahedra add up to $2\pi$ around each edge of $M$, the triangulation is called an \emph{angled triangulation}. Casson proved that every orientable cusped 3--manifold that admits an angled triangulation must also admit a hyperbolic metric, and outlined a possible way to find the hyperbolic metric by studying the volumes of angled tetrahedra --- an idea also developed by Rivin \cite{rivin-volume}.
The power of Casson's approach lies in the fact that the defining equations of an angled triangulation are both linear and local, making angled triangulations relatively easy to find and deform (much easier than to study an actual hyperbolic triangulation, as in 
\cite{neumann-zagier, sakuma-weeks}
or in some aspects of Thurston's seminal approach \cite{thur-notes}). 

Our goal in this paper is to extend this approach to larger and more complicated building blocks. These blocks can be ideal polyhedra instead of tetrahedra, but they may also have non-trivial topology. In general, an \emph{angled block} will be a $3$--manifold whose boundary is subdivided into faces looking locally like the faces of an ideal polyhedron (in a sense to be defined). The edges between adjacent faces carry prescribed dihedral angles. In Section \ref{sec:blocks}, we will describe the precise combinatorial conditions that the dihedral angles must satisfy. These conditions will imply the following 
generalization of a result by Lackenby \cite[Corollary 4.6]{lack-surg}.

\begin{theorem}\label{thm:block-hyperbolic}
Let $(M, \bdy M)$ be an orientable 3--manifold, subdivided into finitely many 
angled blocks in such a way that the dihedral angles at each edge of $M$ sum to $2\pi$. Then $\bdy M$ consists of tori, and 
the interior of $M$ admits a complete hyperbolic metric.
\end{theorem}

One can prove that a particular manifold with boundary is hyperbolic in a spectrum of practical ways, ranging from local to global. In some cases, a combinatorial description of $M$ naturally guides a way to subdivide it into tetrahedra (see, for example, \cite{gf-bundle} or  \cite{weeks:hyp-knots}).
In these cases, angled triangulations are highly useful. On the other extreme, one can study the global topology of $M$ and prove that it contains no essential spheres, disks, tori, or annuli; Thurston's hyperbolization theorem then implies that $M \setminus \bdy M$ is hyperbolic \cite{thur-survey}.
Theorem \ref{thm:block-hyperbolic} provides a medium--range solution 
(still relying on Thurston's theorem)
for situations where $M$ naturally decomposes into pieces that retain some topological complexity. 

We will apply Theorem \ref{thm:block-hyperbolic} to the complements of arborescent links, which are defined in terms of \emph{bracelets}. 
We choose an orientation of $\mathbb{S}^3$, to remain fixed throughout the paper.

\begin{figure}[h]
\psfrag{d0}{$d=0$}
\psfrag{d1}{$d=1$}
\psfrag{d2}{$d=2$}
\psfrag{d3}{$d=3$}
\begin{center}
\includegraphics{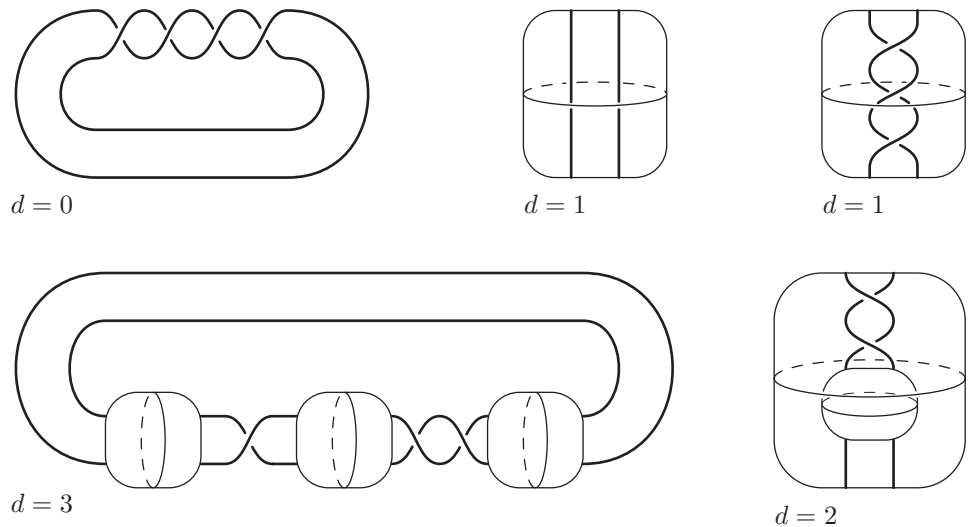}
\end{center}
\vspace{-1ex}
\caption{Examples of $d$--bracelets. The two $1$--bracelets with different numbers of half-twists in their bands are homeomorphic.}
\label{fig:bracelets}
\end{figure}

\begin{define}\label{def:bracelet}
An \emph{unknotted band} $A \subset \mathbb{S}^3$ is an annulus or M\"obius band, whose core curve $C$ is an unknotted circle. Such an $A$ has a natural structure as an $I$--bundle over $C$, and we will refer to the fiber over a point of $C$ as a \emph{crossing segment} of the unknotted band $A$.

Consider the manifold $M_d$ obtained by removing from $\mathbb{S}^3$ the open regular neighborhoods of $d$ disjoint crossing segments of an unknotted band $A$, and let $K_d = \partial A \cap M_d$. Then a \emph{$d$--bracelet} $B_d$ is the pair $(M_d, K_d)$, as in Figure \ref{fig:bracelets}. We say that $d$ is the \emph{degree} of the bracelet.
\end{define}

Note that when $d>0$, $B_d$ is determined up to homeomorphism (of pairs) by the integer $d$. For example, when $d=2$, $B_2$ is homeomorphic to
the pair $(\mathbb{S}^2 \cross I, \{\mbox{4 points}\} \cross I)$.
When $d=1$, $M_1$ is a 3--ball and $K_1$ is a pair of 
simultaneously boundary--parallel arcs; a 1--bracelet $B_1$ is commonly called a \emph{trivial
tangle}. When $d=0$, $B_0$ is determined by the \emph{number of half-twists} in the band: namely, the linking number of $C$ with $\bdy A$.

\begin{figure}%[ht]
\begin{center}
\includegraphics{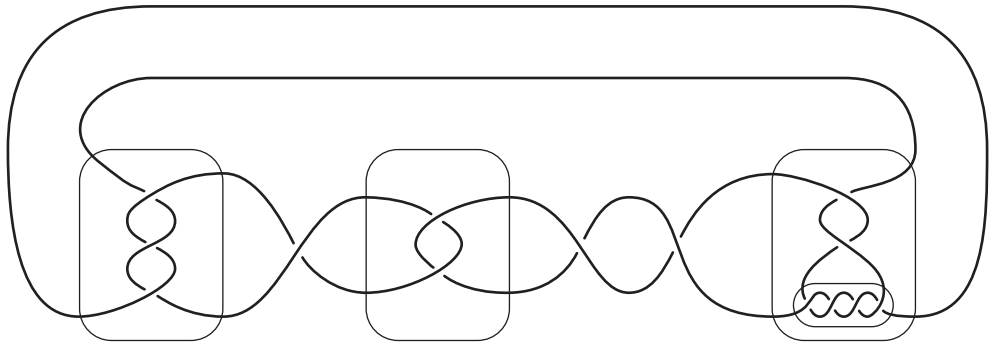}
\caption{A generalized
arborescent knot, obtained by gluing several bracelets.}
\label{fig:arb-link}
\end{center}
\end{figure}

Let $B_{d_1}$ and $B_{d_2}$ be two bracelets with $d_i>0$,
and choose a boundary sphere $S_i$ of each $B_{d_i}$. The $S_i$ have natural orientations
induced by the orientation of $\mathbb{S}^3$,
and we can glue $S_1$ to $S_2$ by any orientation--reversing homeomorphism sending the 
unordered
$4$-tuple of points $S_1\cap K_1$ to the $4$-tuple $S_2\cap K_2$. The union of the $K_{d_i}$ then defines 
a collection of arcs in a larger subset of $\mathbb{S}^3$. More generally, if bracelets $B_{d_1}, \ldots, B_{d_n}$ are glued to form $\mathbb{S}^3$ (some of the $d_i$ being $1$), the arcs in these bracelets combine to form a link $K$ in $\mathbb{S}^3$, as in Figure \ref{fig:arb-link}. 

\begin{define}\label{def:prime}
A link $K \subset \mathbb{S}^3$ is called \emph{prime} if, for every 2--sphere $S$ meeting $K$ in two points, at least one of the two balls cut off by $S$ intersects $K$ in a single boundary--parallel arc. If $K$ is not prime, it is called \emph{composite}. Note that with this convention, every split link (apart from the split link consisting of two unknots) is automatically composite.
\end{define}

\begin{define}\label{def:arborescent}
A knot or link $K = \bigcup_{i=1}^n K_{d_i}$, obtained when several bracelets are glued together to form $\mathbb{S}^3$, is called a \emph{generalized arborescent link}. If, in addition, $K$ is prime, we say that it is an \emph{arborescent link}. \end{define}

The pattern of gluing bracelets to form a link can be represented by a tree $T$, in which a $d$--valent vertex corresponds to a $d$--bracelet and an edge corresponds to a gluing map of two neighboring bracelets. The term \emph{arborescent}, from the Latin word \emph{arbor} (tree), refers to this correspondence. Special cases of arborescent links include \emph{two--bridge links}, which can be constructed by gluing two $1$--bracelets,
and \emph{Montesinos links}, which can be constructed by gluing a single $d$--bracelet to $d$ different $1$--bracelets.
Montesinos links are also known as \emph{star links}, because the corresponding tree is a star.

The tree that represents an arborescent link carries a great deal of geometric and topological information. For example, Gabai has used trees to construct an algorithm that computes the genus of an arborescent link \cite{gabai:arborescent}. Bonahon and Siebenmann have used trees to completely classify arborescent links up to isotopy \cite{bonsieb-monograph}.
One geometric consequence of their work is the following result.

\begin{figure}[ht]
\psfrag{p}{$p$}
\psfrag{q}{$q$}
\psfrag{r}{$r$}
\psfrag{1}{I}
\psfrag{2}{II}
\psfrag{3}{III}
\begin{center}
\includegraphics{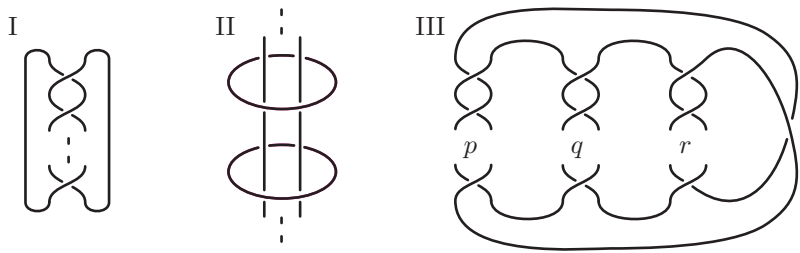}
\end{center}
\caption
{The three exceptional families of non-hyperbolic arbor\-escent links. For family III, $p,q,r\geq2$ and $\frac{1}{p} + \frac{1}{q} + \frac{1}{r} \geq 1$.}
\label{fig:arb-exceptions}
\end{figure}

\begin{theorem}[Bonahon--Siebenmann]\label{thm:main}
The following three families, shown in \linebreak Figure \ref{fig:arb-exceptions}, form a complete list of non-hyperbolic arborescent links:
% DF: hard pagebreak to format the theorem correctly. Remove if we edit the text!
\pagebreak
%%%%%%%%%%%%%
\begin{enumerate}
\item[I.] $K$ is the boundary of a single unknotted band,
\item[II.] $K$ has two isotopic components, each of which bounds a 2--punctured disk
properly embedded in $\mathbb{S}^3\smallsetminus K$,
\item[III.] $K$ or its reflection is the pretzel link $P(p,q,r,\minus 1)$, where $p,q,r\geq2$ and $\frac{1}{p} + \frac{1}{q} + \frac{1}{r} \geq 1$.
\end{enumerate}
Furthermore, an effective algorithm decides whether a given generalized arborescent link $K$ is prime, and whether it lies in one of the exceptional families.
\end{theorem}

Bonahon and Siebenmann's original proof of this theorem made strong
use of the double branched covers of arborescent links. These covers are all \emph{graph manifolds}, obtained by gluing Seifert fibered manifolds along incompressible tori that project to gluing spheres of $d$--bracelets. 
Their results and ideas were heavily quoted, but unfortunately the monograph containing the proof \cite{bonsieb-monograph} has never been finished. One of our primary motivations in this paper was to write down a version of the proof.

In the years since Bonahon and Siebenmann's monograph, several authors have re-proved parts of the theorem. Menasco \cite{menasco-alt} proved that a two--bridge link (more generally, a prime alternating link) is hyperbolic whenever it is not in family I. Oertel  \cite{oertel:star-link} proved that the complement of a Montesinos link contains an incompressible torus if and only if the link is in family III, with $\frac{1}{p} + \frac{1}{q} + \frac{1}{r} = 1$. 
Finally, it follows from Wu's work on Dehn surgery \cite{wu-knot} that all non-Montesinos arborescent \emph{knots} are hyperbolic. 

It is fairly straightforward to check that the links listed in Theorem \ref{thm:main} are indeed non-hyperbolic. For families I and II, Figure \ref{fig:arb-exceptions} reveals an obvious annulus or M\"obius band that forms an obstruction to 
the existence of
a hyperbolic structure. Meanwhile, the pretzel links in family III contain (less obvious) incompressible tori when $\frac{1}{p} + \frac{1}{q} + \frac{1}{r} = 1$ (by Oertel's work \cite{oertel:star-link})
and are Seifert fibered by Sakuma's work \cite{sakuma:star-link} when $\frac{1}{p} + \frac{1}{q} + \frac{1}{r} > 1$ (in fact, such links are torus links unless $(p,q,r)$ is a permutation of $(2,2,n)$).
In particular, all of these well--studied links are known to be prime.
Thus we will focus our attention on proving that all the remaining arborescent links are indeed hyperbolic.

The proof is organized as follows. In Section \ref{sec:blocks}, we will define angled blocks and prove Theorem \ref{thm:block-hyperbolic}. In Section \ref{sec:algorithm}, we will perform a detailed study of how $d$--bracelets can be glued along 4--punctured spheres. This will enable us to simplify the 
bracelet
presentation of any particular link and decide whether it is an exception. In Section \ref{sec:complement}, we will use the bracelet structure to subdivide the link complement into tetrahedra and solid tori. The subdivision will work for all arborescent links except families I and II.
Finally, in Section \ref{sec:angles}, we will assign dihedral angles to edges on the boundary of the tetrahedra and solid tori. For links that are not in family III, these angles will satisfy the criteria of angled blocks, implying by Theorem \ref{thm:block-hyperbolic} that the link complement is hyperbolic.

\smallskip
{\bfseries Acknowledgements:} This project began as the first author's Ph.\ D. thesis under the guidance of Steve Kerckhoff, was nourished by advice from Francis Bonahon, and reached its completion while both authors were visiting Osaka University and enjoying the hospitality of Makoto Sakuma. All three of these mentors deserve our deep
gratitude for their help and encouragement. We would also like to thank Fr\'ed\'eric Paulin for his careful reading and suggestions.

 \section{Angled blocks}\label{sec:blocks}
 
In this section, we develop a theory of angled blocks that provides a practical way of proving that a given manifold is hyperbolic (Theorem \ref{thm:block-hyperbolic}). We lay out the necessary definitions in Section \ref{sec:angled-def}. In Section \ref{sec:normal-surf}, we study the intersections between blocks and surfaces in a manifold, and prove that any surface can be placed into a sufficiently nice \emph{normal form}. The angle structures on the blocks allow us to define a natural measure of complexity for the surfaces, called \emph{combinatorial area}, which behaves like hyperbolic area.
In Section \ref{sec:comb-area}, we will use combinatorial area considerations to show that $M$ cannot contain any essential surfaces of non-negative Euler characteristic, so by Thurston's hyperbolization theorem $M$ must admit a hyperbolic structure.

Our proof of Theorem \ref{thm:block-hyperbolic} follows the same outline as
Casson's proof that manifolds with an angled triangulation are hyperbolic, written down by Lackenby in \cite[Section 4]{lack-surg}. The credit for developing these ideas 
goes mainly to Casson and Lackenby.
 
 \subsection{From polyhedra to blocks}\label{sec:angled-def}
In studying a 3--manifold $M$, it is frequently useful to decompose $M$ into pieces that are not contractible. This idea has been recently studied by other authors: Agol has described a way to cut a manifold into non-contractible \emph{nanotubes} \cite{agol-nanotube}, while Martelli and Petronio have cut a manifold into \emph{bricks} \cite{martelli-petronio}. Rieck and Sedgwick, among others, have investigated how
a solid torus added during Dehn surgery can intersect a Heegaard surface \cite{rieck-sedgwick}. Focusing on the individual pieces of the decomposition, Schlenker has studied manifolds with polyhedral boundary \cite{schlenker}. Our angled blocks fit into this theme.
 
 \begin{define}\label{def:dual-graph}
 Let $S$ be a closed oriented surface,
and let $\Gamma \subset S$ be an embedded graph
each of whose vertices has degree at least $3$.
We say that $\Gamma$ \emph{fills} $S$ if every component of $S \setminus \Gamma$ is an open disk, whose boundary consists of at least $3$ edges of $\Gamma$. Given a graph $\Gamma$ that fills a surface, we can construct a \emph{dual graph} $\Gamma^*
\subset S$, well-defined up to isotopy,
in the following fashion. Every disk of $S \setminus \Gamma$ defines a vertex of $\Gamma^*$. Every edge $e \subset \Gamma$ separates two faces of $S \setminus \Gamma$; we connect the corresponding vertices of $\Gamma^*$ by a dual edge $e^*$.
Finally, $S\smallsetminus \Gamma^*$ is a union of disks, or faces, each corresponding to a vertex of $\Gamma$.

Note that this construction still makes sense if the surface $S$ has several components. In this situation, both $\Gamma$ and $\Gamma^*$ will have as many components as $S$.
\end{define}

\begin{define}\label{def:angled-block}
Let $\pp$ be a compact, oriented, irreducible, atoroidal 3--manifold with boundary. Let $\Gamma$ be a graph that fills $\bdy \pp$, whose edges are $e_1, \ldots, e_n$. To every edge $e_i \subset \Gamma$ we assign an \emph{internal angle} $\alpha_i$ and an \emph{external angle} $\varepsilon_i = \pi - \alpha_i$. By duality, an edge $e^*_i \subset \Gamma^*$ receives the same angle as its dual edge $e_i \subset \Gamma$.

We say that $\pp$ is an \emph{angled block} if this assignment of angles satisfies the following properties:
\begin{enumerate}
\item $0 < \alpha_i < \pi$ for all $i$,
\item $\sum_{\bdy D} \varepsilon_i = 2\pi$ for every face $D$ of $\bdy \pp \setminus \Gamma^*$, and
\item $\sum_{\gamma} \varepsilon_i > 2\pi$ for every simple closed curve $\gamma \subset \Gamma^*$ that bounds a disk in $\pp$ but is not the boundary of a face of $\bdy \pp \setminus \Gamma^*$.
\end{enumerate}
Finally, we remove from $\pp$ all the vertices of $\Gamma$, making them into \emph{ideal vertices}. We  will refer to the edges of $\Gamma$ as the \emph{edges of $\pp$}, and to the faces of $\bdy \pp \setminus \Gamma$ as the \emph{faces of $\pp$}. Removing the vertices of $\Gamma$ makes the faces of $\pp$ into ideal polygons.
\end{define}

Property $(1)$ says that $\pp$ is locally convex at every edge. Property $(2)$ says that the link of every ideal vertex of $\pp$ has the angles of
a convex Euclidean polygon. Property $(3)$ is motivated by the following theorem of Rivin \cite{rivin-idealpoly}:

\begin{theorem}[Rivin] \label{thm:rivin-poly}
Let $\mathcal{P}$ be an angled polyhedron --- that is, a contractible angled block. Then $\mathcal{P}$ can be realized as a convex ideal polyhedron in $\HH^3$ with the prescribed dihedral angles, uniquely up to isometry. Conversely, the dihedral angles of every convex ideal polyhedron in $\HH^3$ satisfy $(1)$--$(3)$.
\end{theorem}

Such characterizations of polyhedra in $\mathbb{H}^3$ by their dihedral angles
were first studied by Andreev \cite{andreev}. We conjecture that an analogous result holds for non-contractible blocks as well:

\begin{conjecture}\label{conj:block-geometry}
Let $\pp$ be an angled block. Then its universal cover $\tilde{\pp}$ can be realized as a (possibly infinite) ideal polyhedron in $\HH^3$, with dihedral angles specified by $\pp$, uniquely up to isometry.
\end{conjecture}

When $\pp$ is contractible, this conjecture is exactly
Rivin's theorem. Schlenker \cite[Theorem 8.15]{schlenker} has treated the case where $\pp$ has incompressible boundary.
Finally, when all angles are of the form $\pi/n$ with $n\geq 2$, the conjecture follows by a doubling argument from the hyperbolization theorem for orbifolds \cite{blp, chk}.
We also note that the converse statement (that the ideal polyhedron $\tilde{\pp}$ must satisfy $(1)$--$(3)$) is a fairly straightforward consequence of the Gauss--Bonnet theorem.

Our primary interest is in the manifolds that one may construct by gluing together angled blocks. To build a manifold with boundary, we first truncate all the ideal vertices of the blocks. As a result, a block $\pp$ has two kinds of faces: \emph{interior faces} that are truncated copies of the original faces, and \emph{boundary faces} that come from the truncated vertices. Similarly, $\pp$ has two kinds of edges: \emph{interior edges} that are truncated edges of $\Gamma$, and \emph{boundary edges} along the boundary faces. We note that a truncated block is a special case of a differentiable \emph{manifold with corners} (modeled over $\mathbb{R}_+^3$: see \cite{douady} for a general definition).

\begin{define}\label{def:angled-decomposition}
Let $(M, \bdy M)$ be a compact $3$--manifold with boundary. An {\it angled
decomposition of $M$} is a subdivision of 
$M$ into truncated angled blocks, glued along their interior faces, such that $\sum \alpha_i = 2 \pi$ around each interior edge of $M$. The boundary faces of the blocks fit together to tile $\bdy M$.
\end{define}

Theorem \ref{thm:block-hyperbolic} says that the interior of
every orientable manifold with an angled decomposition must admit a hyperbolic structure. However, this is purely an existence result. An angled decomposition of a manifold is considerably weaker and more general than a hyperbolic structure, for two reasons. First, we do not know whether the blocks are actually geometric pieces --- this is the content of Conjecture \ref{conj:block-geometry}. Second, even when the blocks are known to be geometric, a geometrically consistent gluing must respect more than the dihedral angles. To obtain a complete hyperbolic structure, the truncated vertices of the blocks must fit together to tile a horospherical torus, meaning that these Euclidean polygons must have consistent sidelengths as well as consistent angles.

There is an interesting contrast between the rigidity of a hyperbolic structure and the flexibility of angle structures. By Definitions \ref{def:angled-block} and \ref{def:angled-decomposition}, an angle structure
on a block decomposition is a solution to a system of linear equations and (strict) linear inequalities. The solution set to this system, if non-empty, is an open convex polytope, so for every angled decomposition there is a continuum of deformations.

In fact, geometric angled blocks --- for example, angled polyhedra --- can serve as a stepping stone on the way to finding a complete hyperbolic structure. Every angled polyhedron has a well-defined volume determined by its dihedral angles, by Theorem \ref{thm:rivin-poly}.
If the volume of an angled decomposition is critical in the polytope of deformations, we can exploit 
Schl\"afli's formula as in Rivin's theorem \cite{rivin-volume} and show that the polyhedra glue up to give a hyperbolic metric: this is carried out for some examples in \cite{gf-bundle} (where all blocks are tetrahedra).
However, depending on the combinatorics of the decomposition, a critical point may or may not occur. In fact, numerical experiments show that some of the decompositions that we will define for arborescent link complements in Section \ref{sec:complement} admit angle structures, but have no critical point.

\subsection{Normal surface theory in angled blocks}\label{sec:normal-surf}

To prove Theorem \ref{thm:block-hyperbolic}, we study the intersections between angled blocks and (smooth) essential surfaces.

\begin{define}\label{def:essential}
A surface $(S, \bdy S) \subset (M, \bdy M)$ is called \emph{essential} if $S$ is incompressible, boundary--incompressible, and not boundary--parallel, or if $S$ is a sphere that does not bound a ball.
\end{define}

Our goal is to move any essential surface into a form where its intersections with the individual blocks are particularly nice:

\begin{define}\label{def:normal-surf}
Let $\pp$ be a truncated block, and let $(S, \bdy S) \subset (\pp, \bdy \pp)$ be a surface. We say that $S$ is \emph{normal} if it satisfies the following properties:
\begin{enumerate}
\item every closed component of $S$ is essential in $\pp$,
\item $S$ and $\partial S$ are transverse to all faces and edges of $\mathcal{P}$,
\item no component of $\bdy S$ lies entirely in a face of $\bdy \pp$,
\item no arc of $\bdy S$ in a face of $\pp$ runs from an edge of $\pp$ back to the same edge,
\item no arc of $\bdy S$ in an interior face of $\pp$ runs from a boundary edge to an adjacent interior edge.
\end{enumerate}
Given a decomposition of $M$ into blocks, a surface $(S, \bdy S) \subset (M, \bdy M)$ is called \emph{normal} if for every block $\pp$, the intersection $S \cap \pp$ is a normal surface
in $\mathcal{P}$.
\end{define}

\begin{theorem}\label{lemma:normalize}
Let $(M, \bdy M)$ be a manifold with an angled block decomposition.  
\begin{itemize}
\item[(a)]  If $M$ is reducible, then $M$ contains a normal 2--sphere.
\item[(b)] If $M$ is irreducible and $\bdy M$ is compressible, then $M$ contains a normal disk.
\item[(c)] If $M$ is irreducible and $\bdy M$ is incompressible, then any essential surface can be moved by isotopy into normal form. 
\end{itemize}
\end{theorem}

\begin{proof}
The following argument is the standard procedure for placing surfaces in normal form with respect to a triangulation or a polyhedral decomposition \cite{haken-normal}. As long as all faces of all blocks are disks, the topology of the blocks never becomes an issue. We will handle part (c) first, followed by (b) and (a). 

\medskip

%%%%%% Incompressible, irreducible

For (c), assume that $M$ is irreducible and $\bdy M$ is incompressible. Let $(S, \bdy S)$  be an essential surface in $(M, \bdy M)$. To move $S$ into normal form, we need to check the conditions of Definition \ref{def:normal-surf}. Since $S$ is essential in $M$, it automatically satisfies $(1)$. Furthermore, a small isotopy of $S$
ensures the transversality conditions of $(2)$.

Consider the intersections between $S$ and the open faces of the blocks, and let $\gamma$ be one component of intersection. Note that by Definition \ref{def:dual-graph}, the face $F$ containing $\gamma$ is contractible. We want to make sure that $\gamma$ satisfies $(3)$, $(4)$, and $(5)$.

\begin{enumerate}
\setcounter{enumi}{2}

\item Suppose that $\gamma$ is a closed curve, violating $(3)$. Without loss of generality, we may assume that $\gamma$ is innermost on the face $F$. Then $\gamma$ bounds a disk $D \subset F$, whose interior is disjoint from $S$. But since $S$ is incompressible, $\gamma$ also bounds a disk $D' \subset S$. Furthermore, since we have assumed that $M$ is irreducible, the sphere $D \cup_\gamma D'$ must bound a ball. Thus we may isotope $S$ through this ball, moving $D'$ past $D$. This isotopy removes the curve $\gamma$ from the intersection between $S$ and $F$. 

\smallskip

\begin{figure}[h]
\psfrag{g}{$\gamma$}
\psfrag{e}{$e$}
\psfrag{D}{$D$}
\psfrag{dp}{$D'$}
\psfrag{S}{$S$}
\psfrag{bs}{$\bdy S$}
\psfrag{bm}{$\bdy M$}
\begin{center}
\includegraphics{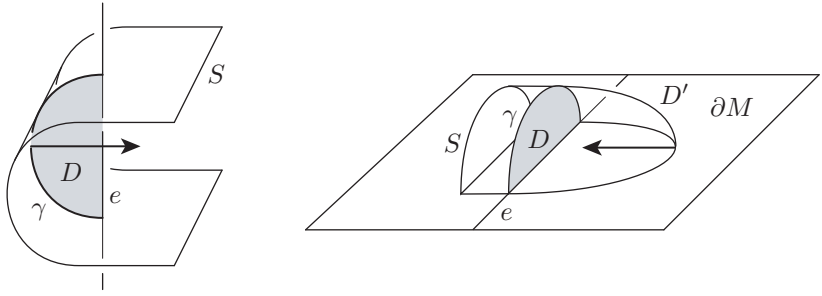}
\end{center}
\caption{When a surface violates condition $(4)$ of normality, then an isotopy in the direction of the arrow removes intersections between $S$ and the faces of $M$.}
\label{fig:normalize}
\end{figure}

\item Suppose that $\gamma$ runs from an edge $e$ back to $e$, violating $(4)$.
Then $\gamma$ and $e$ co-bound a disk $D \subset F$, and we can assume $\gamma$ is innermost (i.e. $S$ does not meet $D$ again).
If $e$ is an interior edge, we can use this disk $D$ to guide an isotopy of $S$ past the edge $e$, as in the left panel of Figure \ref{fig:normalize}. This isotopy removes $\gamma$ from the intersection between $S$ and $F$
% FG: add
(some intersection components between $S$ and the interiors of faces other than $F$ may merge, but their total number 
% DF: clarify
% never increases).
always decreases).
% FG: "always decreases" is nicer to the ear, but the number of intersection components with faces "other than F" could stay constant. Decrementation comes from F. Oh well...
%

If $\gamma$ lies in a boundary face, then the situation is very similar to the previous paragraph. This time, the disk $D$ guides an isotopy of $S$ along $\bdy M$, simplifying the intersection between $S$ and the faces of the blocks.

Finally, if $e$ is a boundary edge and $F$ is an interior face, then $D$ is a boundary compression disk for $S$. Since $S$ is boundary--incompressible, $\gamma$ must also cut off a disk $D' \subset S$, as in the right panel of Figure \ref{fig:normalize}. Since $M$ is irreducible and $\bdy M$ is incompressible, it follows that the disk $D \cup_\gamma D'$ is boundary--parallel:
$D\cup D'\cup \Delta$ bounds a ball $B$, for some disk $\Delta \subset \partial M$. We must ask on which side of $D\cup_{\gamma} D'$ the ball $B$ lies: if a neighborhood of the arc $\gamma$ in the surface $S$ meets the interior of $B$, then $S$ is a disk of $B$ and is boundary--parallel (recall that $S$ does not meet $D$ again, because $\gamma$ is innermost among the arcs running from $e$ back to $e$). So $S$ does not meet the interior of $B$. In particular, $D'$ is isotopic to $D$ by an isotopy sweeping out $B$ and missing $S\smallsetminus D'$. This defines an isotopy of $S$ which can be extended slightly to move $D'$ past $D$, thus removing the curve $\gamma$ from the intersection between $S$ and $F$.

\smallskip

\begin{figure}[h]
\psfrag{g}{$\gamma$}
\psfrag{e}{$e$}
\psfrag{D}{$D$}
\psfrag{S}{$S$}
\psfrag{bs}{$\bdy S$}
\psfrag{bm}{$\bdy M$}
\begin{center}
\includegraphics{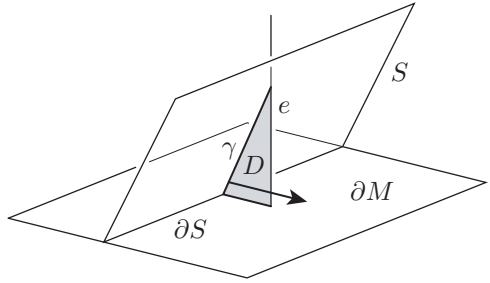}
\end{center}
\caption{When a surface violates condition $(5)$ of normality, a $\bdy M$--preserving isotopy of $S$ along the disk $D$, in the direction of the arrow, removes intersections between $S$ and the faces of $M$.}
\label{fig:normalize2}
\end{figure}

\item Suppose that $\gamma$ runs from a boundary edge to an adjacent interior edge, violating $(5)$. Then $\gamma$ once again cuts off a disk $D$. By isotoping $S$ along this disk, as in Figure \ref{fig:normalize2}, we remove $\gamma$ from the intersection.
\end{enumerate}

It is immediate to check that
each of the last three moves reduces the number of
components of $S\cap Z$, where $Z$ is the union of the interiors of the faces of $M$. Thus, after a finite number of isotopy moves, $S$ becomes normal.

\medskip

%%% Compressible, irreducible

For part (b), assume that $M$ is irreducible and $\bdy M$ is compressible. Let $S$ be an essential disk in $M$; under our assumptions, $S$ must be a compression disk for $\bdy M$. To move $S$ into normal form, we follow a very similar procedure to the one in part (c). In particular, condition $(1)$ of Definition \ref{def:normal-surf} is vacuous because $S$ has no closed components. Furthermore, a small isotopy of $S$ ensures the transversality conditions of $(2)$. Focusing our attention on conditions $(3)-(5)$, let $\gamma$ be one component of intersection between $S$ and a face $F$ of a block.

If $\gamma$ is a simple closed curve, violating $(3)$, the argument is exactly the same as above. We find that $\gamma$ bounds a disk $D \subset F$ and an isotopic disk $D'  \subset S$, because $S$ is incompressible and $M$ is irreducible. Thus we may isotope $S$ past $D$.

If $\gamma$ runs from an edge $e$ back to $e$, violating $(4)$, the argument is mostly the same as above. If $e$ is an interior edge, or $\gamma$ lies in a boundary face, then the exact isotopies described in part (a)--$(4)$ 
will guide $S$ past $e$. If $e$ is a boundary edge and $\gamma$ lies in an interior face $F$, then $\gamma$ and $e$ co-bound a disk $D\subset F$;
up to replacing $\gamma$ with an outermost arc of $D\cap S$ on $F$, we may assume $D\cap S=\gamma$ so that the disk $D$
realizes a boundary compression of $S$. The situation is similar to the right panel of Figure \ref{fig:normalize}, except now $\gamma$ splits $S$ into disks $D_1$ and $D_2$ (since $S$ itself is a disk). At least one $D_i \cup_\gamma D$ must be essential in $M$, because if they were both boundary--parallel, $S$ would be boundary--parallel also. 
If by $S$ we now denote this essential disk $D_i \cup_\gamma D$, then $S$
can be pushed away from the face $F$.

Finally, if $\gamma$ runs from a boundary edge to an adjacent interior edge $e$, violating $(5)$, an isotopy of $S$ as in Figure \ref{fig:normalize2} will remove $\gamma$ from the intersections between $S$ and the faces of the blocks. Each of the last three moves simplifies the intersections between $S$ and the faces, so a repeated application will place $S$ in normal form.
\medskip

%% Reducible

For part (a), assume that the manifold $M$ is reducible, and let $S \subset M$ be a sphere that doesn't bound a ball. We will move $S$ into normal form by checking the conditions of Definition \ref{def:normal-surf}. Note that by Definition \ref{def:angled-block}, an essential sphere can never be contained in a single block, so condition $(1)$ is vacuous. A small isotopy of $S$
ensures the transversality conditions of $(2)$. Note as well that condition $(5)$ is vacuous, because $S$ is closed. To satisfy conditions $(3)$ and $(4)$, let $\gamma$ be one arc of intersection between $S$ and a face $F$ of a block. 

If $\gamma$ is a simple closed curve, violating $(3)$, we may assume as before that $\gamma$ is innermost in $F$. Thus $\gamma$ bounds a disk $D \subset F$ whose interior is disjoint from $S$. Because $S$ is a sphere, we may write $S=D_1 \cup_{\gamma} D_2$ for disks $D_1$ and $D_2$. Suppose that each $D_i \cup_{\gamma} D$ bounds a ball $B_i$. Because the boundaries of  $B_1$ and $B_2$ intersect exactly along a single disk $D$, either one ball contains the other or they have disjoint interiors. In either scenario, it follows that $S = D_1 \cup D_2$ must bound a ball --- a contradiction. Thus, since at least one $D_i \cup_{\gamma} D$ must fail to bound a ball, we can replace $D$ by one of the $D_i$. The resulting sphere, which we continue to call $S$, can be pushed away from the face $F$.

If $\gamma$ runs from an edge $e$ back to $e$, violating $(4)$, then $\gamma$ and $e$ co-bound a disk $D$. As before, we can use $D$ to guide an isotopy of $S$ past $e$. (See Figure \ref{fig:normalize}, left.) Note that since $S$ is closed, $\gamma$ must be an interior edge.

By repeating these moves, we eventually obtain a sphere in normal form.
\end{proof}

\subsection{Combinatorial area}\label{sec:comb-area}

So far, we have not used the dihedral angles of the blocks. Their use comes in estimating the complexity of normal surfaces.

\begin{define}\label{def:comb-area}
Let $\pp$ be an angled block, and denote by $\varepsilon_{\delta}$ the exterior dihedral angle at the edge $\delta$. Truncate the ideal vertices of $\pp$, and label every boundary edge $\delta$ with a dihedral angle of $\varepsilon_{\delta}=\frac{\pi}{2}$. Let $S$ be a normal surface in $\pp$, and let $\delta_1,\dots, \delta_n$ be the edges of the truncated block $\mathcal{P}$ met by $\partial S$ (each edge may be counted several times). We define the \emph{combinatorial area} of $S$ to be
$$a(S) = \sum_{i=1}^n \varepsilon_{\delta_i} - 2 \pi \chi(S).$$
For the sake of brevity, we will refer to the above sum of dihedral angles 
% FG: insert to make clear the \chi is not counted
($\sum_{i=1}^n \varepsilon_{\delta_i}$)
as $\sum_{\partial S} \varepsilon_i$.
\end{define}

Note that by the Gauss--Bonnet theorem, the right-hand side is just the area of a hyperbolic surface with piecewise geodesic boundary, with 
exterior angles $\varepsilon_i$ 
along the boundary and Euler characteristic $\chi(S)$.

\begin{figure}[h]
\begin{center}
\psfrag{(a)}{(a)}
\psfrag{(b)}{(b)}
\includegraphics{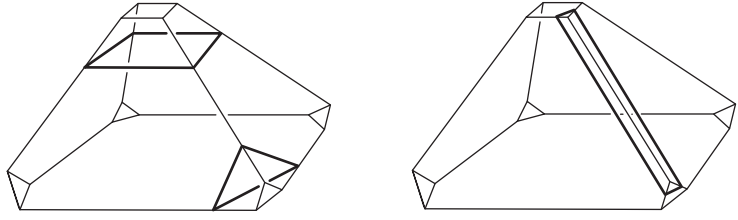}
\caption{In any angled block, \emph{vertex links} (left) and \emph{boundary bigons} (right) are the only connected normal surfaces of area $0$.}
\label{fig:zero-area}
\end{center}
\end{figure}

%%%%%%%%%%%%%%%%%%%%%%%%%%
\begin{lemma}\label{lemma:pos-area}
Let $S$ be a normal surface in a truncated angled block $\pp$. Then $a(S) \geq 0$. Furthermore, if $a(S)=0$, then every component of $S$ is a vertex link 
(boundary of a regular neighborhood of a boundary face)
or a boundary bigon
(boundary of a regular neighborhood of an interior edge),
as in Figure \ref{fig:zero-area}.
\end{lemma}

\begin{proof}
Because combinatorial area is additive over multiple components of $S$, it suffices to consider the case when $S$ is connected. Furthermore, when $\chi(S)<0$, $a(S)>0$, so it suffices to consider the case when $\chi(S) \geq 0$. By Definition \ref{def:angled-block}, $\pp$ is irreducible and atoroidal, so $S$ cannot be a sphere or torus. If $S$ is an annulus, $a(S) = \sum_{\bdy S} \varepsilon_i >0$, because $\bdy S$ must intersect some edges and the dihedral angle on each edge is positive. Thus the only 
remaining case is when $S$ is a disk. 

For the rest of the proof, let $D \subset \pp$ be a normal disk. We consider three cases, conditioned on $n$, the number of intersections between $D$ and the boundary faces.

\smallskip
\underline{\emph{Case 0: $n=0$.}}  Recall, from Definitions \ref{def:dual-graph} and \ref{def:angled-block}, that every interior face of $\pp$ corresponds to a complementary region of the graph $\Gamma$ and to a vertex of the dual graph $\Gamma^*$. Thus $\bdy D$ defines a closed path $\gamma$ through the edges of $\Gamma^*$; this is a non-backtracking path because no arc of $\bdy D$ runs from an edge back to itself. The path $\gamma$ may pass through an edge multiple times, but it contains a simple closed curve in $\Gamma^*$. Thus, by Definition \ref{def:angled-block}, $\sum_{\bdy D} \varepsilon_i \geq 2 \pi$. Equality can happen only when $\bdy D$ encircles an ideal vertex, in other words when $D$ is a vertex link.

\begin{figure}[h]
\begin{center}
\psfrag{a1}{$a_1$}
\psfrag{a2}{$a_2$}
\psfrag{d}{$D$}
\psfrag{bp}{$D'$}
\psfrag{bf}{boundary face}
\includegraphics{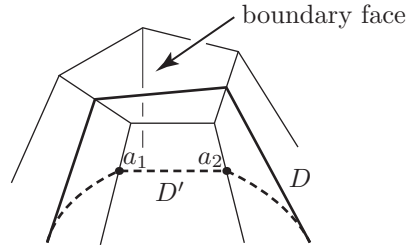}
\caption{We may isotope $\bdy D$ off a boundary face of $\pp$, producing a normal disk $D'$ with $a(D') \leq a(D)$.}
\label{fig:isotope-off-bdy}
\end{center}
\end{figure}

\smallskip
\underline{\emph{Case 1: $n=1$.}}  The two boundary edges crossed by $\bdy D$ contribute $\pi$ to the external angle sum of $\bdy D$. Thus we may isotope $\bdy D$ off the boundary face without increasing the angle sum, since by Definition \ref{def:angled-block} the interior edges meeting this face have a total angle of $2\pi$. Let $D'$ be the resulting disk, and $a_1,\dots,a_k$ be the intersection points, numbered consecutively, of $\partial D'$ with interior edges of the block near the old boundary face (in Figure \ref{fig:isotope-off-bdy}, $k=2$).

We claim that $D'$ is normal, and is not a vertex link. 
Since \linebreak $n=1$,
the only way that $D'$ can fail Definition \ref{def:normal-surf} is if an arc of $\bdy D'$ violates condition $(4)$ and runs from an interior edge $e$ back to itself. This cannot happen between $a_i$ and $a_{i+1}$, otherwise the block would have a monogon face, in contradiction with Definitions \ref{def:dual-graph} -- \ref{def:angled-block}. So condition $(4)$ is violated by an arc starting from $a_1$ in the direction opposite $a_2$ to end on the interior edge $e$ (or by an analogous arc from $a_k$). But then the corresponding arc of $D$ must connect $e$ to an adjacent boundary edge, contradicting condition $(5)$. Similarly, the only way to create a vertex link by pulling an arc of $D$ off a boundary face is if all of $D$ is parallel to that boundary face --- but then $D$ once again violates condition $(5)$. Thus, by Case 0, $a(D) \geq a(D') > 0$.

\begin{figure}[h]
\begin{center}
\psfrag{a1}{$a_1$}
\psfrag{a2}{$a_2$}
\psfrag{ak}{$a_k$}
\psfrag{b1}{$b_1$}
\psfrag{bl}{$b_l$}
\psfrag{dd}{$\partial D$}
\psfrag{dD}{$\partial D'$}
\psfrag{D}{$D'$}
\psfrag{F}{$F'$}
\psfrag{do}{$\dots$}
\includegraphics{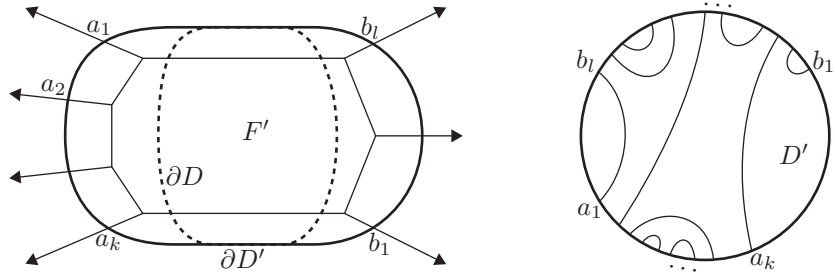}
\caption{When $n=2$ and $a(D)=0$, we have a contradiction for $D'$ normal (left), as well as for $D'$ non-normal (right), unless $k=l=1$.}
\label{fig:find-bigon}
\end{center}
\end{figure}

\smallskip
\underline{\emph{Case 2: $n \geq 2$.}}  Since $\bdy D$ crosses at least 4 boundary edges, $a(D) \geq 0$, with equality only if  $n=2$ and $\bdy D$ is disjoint from the interior edges. We restrict our attention to this case, and claim that $D$ is a boundary bigon.

Push $\partial D$ off the two boundary faces $F_1$ and $F_2$, in a way that minimizes the angle sum of the new disk $D'$. Denote by $a_1, \dots, a_k$ (resp. $b_1,\dots, b_l$) the points where $\partial D'$ crosses interior edges near $F_1$ (resp. $F_2$). Orient the edges containing the $a_i$ and $b_j$ away from the faces $F_1$ and $F_2$. If $A$ (resp. $B$) is the sum of the angles of $D'$ at the $a_i$ (resp. $b_i$), then $A,B\leq\pi$.

Suppose $D'$ is normal. Since we know $a(D') \leq a(D) = 0$, it follows by Case 0 that $D'$ must be the vertex link associated to a boundary face $F'$. Moreover, we have $A=B=\pi$, hence $k,l\geq 2$. Let us isotope $D'$ into $\partial \mathcal{P}$ while keeping its boundary fixed, so that after the isotopy, $D'$ contains the boundary face $F'$ as well as initial segments of all interior edges starting at $F'$: these initial segments end at $a_1, \dots, a_k, b_1, \dots, b_l$, in that cyclic order around $F'$. Suppose the orientations on the interior edges through the $a_i$ are inward for $D'$ (in particular, this will happen whenever $F_1 \neq F'$). Then, since $k\geq 2$, it follows that $\partial \mathcal{P}$ contains an ideal bigon, which is impossible. Therefore the orientations point outward, which implies notably $F_1=F'$. Similarly, $F_2=F'$. As a result, the boundary of the original disk $D$ violated condition $(4)$, e.g. at the boundary edge situated between $a_k$ and $b_1$ (Figure \ref{fig:find-bigon}, left). Contradiction.

Therefore $D'$ is not normal: define $\Delta=D'$ (we are going to modify $\Delta$, but not $D'$). Then the loop $\partial \Delta$ must violate $(4)$, running in a U--turn from an interior edge $e$ back to $e$: we can isotope the disk $\Delta$ so as to erase this U--turn. The angle sum of $\Delta$ decreases to a value less than $2\pi$, so $\Delta$ is even less normal now (by Case 0). If $\partial \Delta$ still crosses any (interior) edges, we can repeat the operation, until $\partial \Delta$ violates $(3)$, and $\Delta$ can be isotoped into an interior face (recall the block $\mathcal{P}$ is irreducible). Therefore $D'$ can be isotoped, with fixed boundary, to a disk in the union of all (open) interior faces and interior edges. Interior edges must connect across $D'$ the points $a_1, \dots, a_k, b_1, \dots, b_l$ of $\partial D'$, which are still cyclically ordered (Figure \ref{fig:find-bigon}, right.) If some edge goes from $a_i$ to $a_j$ (where $i<j$) then there must be an edge from $a_s$ to $a_{s+1}$ for some $i\leq s < j$, and therefore $\partial \mathcal{P}$ contains an ideal monogon: impossible. So every edge across $D'$ runs from an $a_i$ to a $b_j$, in fact to $b_{k+1-i}$ (and we have $k=l$). If $k\geq 2$, then $\partial \mathcal{P}$ contains an ideal bigon. Therefore $k=1$, so $D'$ is traversed by a single edge $e$, and the original disk $D$ was the boundary bigon associated to $e$.
\end{proof}
%%%%%%%%%%%%%%%%%%%%%%%%%%

For an essential surface $(S, \bdy S) \subset (M, \bdy M)$, we can define the combinatorial area $a(S)$ by adding up
the areas of its intersections with the blocks. This definition of combinatorial area was designed to satisfy a Gauss--Bonnet relationship.

\begin{prop}\label{prop:gauss-bonnet}
Let $(S, \bdy S) \subset (M, \bdy M)$ be a surface in normal form. Then
$$a(S) = -2\pi \chi(S).$$
\end{prop}

\begin{proof}
Consider the decomposition of $S$ into $S_1, \ldots, S_n$, namely its components of intersection with the various blocks. Let $S'=S_{i_1}\cup \dots \cup S_{i_k}$ be a union of some $S_i$ glued along some (not necessarily all) of their edges: $S'$ is a manifold with polygonal boundary. Define the interior angle of $S'$ at a boundary vertex to be the sum of the interior angles of the adjacent $S_{i_{\alpha}}$, and the exterior angle as the complement to $\pi$ of the interior angle. It is enough to prove that 
\begin{equation} \sum_{\alpha=1}^{k}a(S_{i_{\alpha}}) \: =\: \sum_{\partial S'}\varepsilon_i -2\pi\chi(S')~, \label{eq:area-additive} \end{equation} 
where the $\varepsilon_i$ are the exterior angles of $S'$: the result will follow by taking $S'=S$ (the union of all $S_i$ glued along all their edges), because all $\varepsilon_i$ are then equal to $\pi-(\frac{\pi}{2}+\frac{\pi}{2})=0$. 
Since $M$ is orientable, up to replacing $S$ with the boundary of its regular neighborhood, we can restrict to the case where $S$ is orientable.

We prove (\ref{eq:area-additive}) by induction on the number of gluing edges,
where the set of involved components $S_{i_{\alpha}}$ is chosen once and for all.
When no edges are glued, (\ref{eq:area-additive}) follows from Definition \ref{def:comb-area}. It remains to check that the right hand side of (\ref{eq:area-additive}) is unchanged when two edges are glued together. In what follows, 
\begin{itemize}
\item $\nu$ is the number of boundary vertices of $S'$,
\item $\theta$ is the sum of all \emph{interior} angles along $\partial S'$, and 
\item $\chi$ is the Euler characteristic of $S'$.
\end{itemize}
Thus the right hand side of (\ref{eq:area-additive}) is $\nu\pi-\theta-2\pi\chi$.

If we glue edges $ab$ and $cd$, where $a,b,c,d$ are distinct vertices of $S'$, then $\theta$ is unchanged, but $\nu$ goes down by $2$ and $\chi$ goes down by $1$: (\ref{eq:area-additive}) is preserved.

If we glue edges $ab$ and $bc$ by identifying $a$ and $c$, where $a,b,c$ are distinct vertices, then $\theta$ goes down by $2\pi$, because $b$ becomes an interior vertex, while $\nu$ goes down by $2$ and $\chi$ is unchanged: (\ref{eq:area-additive}) is preserved.

If we glue two different edges of the form $ab$, closing off a bigon boundary component of $S'$, then $\theta$ goes down by $4\pi$ because both $a$ and $b$ become interior vertices. Since $\nu$ goes down by $2$, and $\chi$ goes up by $1$, again (\ref{eq:area-additive}) is preserved.

If we glue an edge $ab$ to a monogon boundary component $cc$, where $a,b,c$ are distinct vertices, then $\theta$ is unchanged, while $\nu$ goes down by $2$ and $\chi$ goes down by $1$. If we glue two boundary monogons $aa$ and $bb$ together (where $a\neq b$), then $\theta$ goes down by $2\pi$, while $\nu$ goes down by $2$ and $\chi(S')$ is unchanged. In all cases, (\ref{eq:area-additive}) is preserved.
\end{proof}

We are now ready to complete the proof of Theorem \ref{thm:block-hyperbolic}.

\medskip
\noindent{\bf Theorem \ref{thm:block-hyperbolic}.}
\emph{
Let $(M, \bdy M)$ be an orientable $3$--manifold with an angled decomposition. Then $\bdy M$ consists of tori, and the interior of $M$ is hyperbolic.
}
\smallskip

\begin{proof}
Each component of $\bdy M$ is tiled by boundary faces of the
blocks. Just inside each boundary face, a block has a normal
disk of area $0$. These vertex links glue up to form a closed,
boundary--parallel normal surface $S$ of area $0$.  By Proposition
\ref{prop:gauss-bonnet}, $\chi(S) = 0$, and since $M$ is orientable, the
boundary--parallel
surface $S$ must be a torus. Thus $\bdy M$ consists of tori.

By Thurston's hyperbolization theorem \cite{thur-survey}, 
the interior of $M$ 
carries a complete, finite--volume hyperbolic metric
if and only if $M$ contains no essential spheres, disks, annuli, or tori. 
By Theorem \ref{lemma:normalize}, if $M$ has such an essential surface, then it has one in normal form. A normal sphere or disk has positive Euler characteristic, hence negative area.  Thus it cannot occur. 

A normal torus $T \subset M$ has area $0$ and thus, by Lemma
\ref{lemma:pos-area}, must be composed of normal disks of area $0$. Since
$T$ has no boundary, these must all be vertex links, which glue up to
form a boundary--parallel torus. Similarly, a normal annulus $A \subset
M$ must be composed entirely of bigons, since a bigon cannot be glued
to a vertex link. But a chain of bigons forms a tube around an edge of
$M$, which is certainly not essential. Thus we can conclude that $M$
is hyperbolic.
\end{proof}

\section{A simplification algorithm for arborescent links}\label{sec:algorithm}

Recall, from the introduction, that a generalized arborescent link is constructed by gluing together a number of $d$--bracelets. In this section, we describe an algorithm that takes a particular link and simplifies its bracelet presentation into a reduced form. This algorithm, 
directly inspired by Bonahon and Siebenmann's work \cite{bonsieb-monograph},
has several uses. Firstly, if a given generalized arborescent link is composite, the algorithm will decompose it into its prime arborescent pieces. Secondly, the simplified bracelet description will allow us to rapidly identify the non-hyperbolic arborescent links listed in Theorem \ref{thm:main}. In particular, the algorithm recognizes the unknot from among the family of generalized arborescent links. Finally, the simplified bracelet form of an arborescent link turns out to be the right description for \linebreak the block decomposition of the link complement that we 
undertake in Section \ref{sec:complement}.

\subsection{Slopes on a Conway sphere}
Whenever two bracelets are glued together, they are joined along a 2--sphere that intersects the link $K$ in 4 points. This type of sphere, called a \emph{Conway sphere}, defines a 4--punctured sphere in the link complement. Our simplification algorithm is guided by the way in which gluing maps act on arcs in Conway spheres. 

\begin{define}\label{def:sphere-slope}
Let $S$ be a 4--punctured sphere. An \emph{arc pair} $\gamma
\subset S$ consists of two disjoint, properly embedded arcs $\gamma_1$
and $\gamma_2$, such that $\gamma_1$ connects two punctures of $S$ and
$\gamma_2$ connects the remaining two punctures of $S$. A {\it slope} on
$S$ is an isotopy class of arc pairs,
and is determined by any one of the two arcs.
\end{define}

\begin{figure}[ht]
\begin{center}
\includegraphics{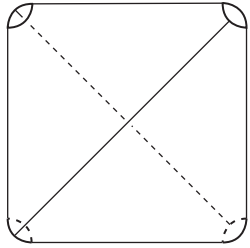}
\caption{Arcs of slope $0$, $1$, and $\infty$ give an ideal
triangulation of a $4$-punctured sphere $S$.}
\label{fig:sphere-slopes}
\end{center}
\end{figure}

To visualize slopes, it helps to picture $S$ as a pillowcase
in $\mathbb{R}^3$
surrounding the unit square of
 $\RR^2$, with punctures at the
corners. (See Figure \ref{fig:sphere-slopes}.) Any arc pair on the pillowcase can then be straightened so
that its intersections with the front of the pillow have a
well-defined Euclidean slope. A \emph{marking} of $S$ (that is, a fixed homeomorphism between $S$ and the pillowcase of Figure \ref{fig:sphere-slopes}) induces a bijection between slopes on $S$ and elements of $\qbar = \QQ \cup \{\infty\}$. 

\begin{figure}[h]
\begin{center}
\includegraphics[width=2.5in]{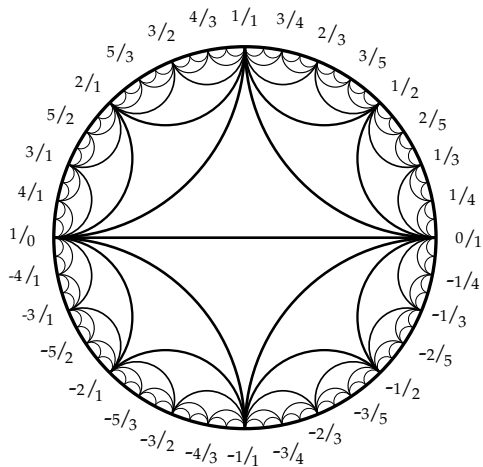}
\caption{The Farey complex $\ff$ of a 4--punctured sphere (graphic by Allen Hatcher).}
\label{fig:farey-graph}
\end{center}
\end{figure}

Slopes on 4--punctured spheres can be neatly represented by the {\it Farey
complex} $\ff$, shown in Figure \ref{fig:farey-graph}. Vertices of $\ff$
correspond to slopes
(arc pairs),
edges of $\ff$ to disjoint slopes, and triangles
to triples of disjoint slopes. Observe that a choice of three disjoint arc pairs of 
different slopes gives an ideal triangulation of $S$. Figure \ref{fig:farey-graph} also illustrates that 
up to a homeomorphism, 
$\ff$ can be identified with the Poincar\'e disk and endowed with a hyperbolic metric, making the triangles of $\ff$ into straight ideal triangles.
The dual of $\ff$ is an infinite trivalent planar tree.

\begin{define}\label{def:preferred-slope}
Let $S$ be a boundary sphere of a $d$--bracelet. This Conway sphere will be assigned a \emph{preferred slope}, as follows. When $d>1$, pick a crossing segment on each side of $S$ (see Definition \ref{def:bracelet}). If we isotope these two segments into $S$, we get an arc pair whose slope is the preferred slope of $S$. When $d=1$, $K_1$ consists of two arcs that can be isotoped into $S$; their slope is then the preferred slope of $S$. Note that the two definitions (for $d=1$ and $d>1$) are truly different. The case $d=0$ is empty (a $0$--bracelet has no boundary spheres).
\end{define} 

\subsection{The algorithm}

We will perform the following sequence of steps to simplify the bracelet presentation of a generalized arborescent link. 

\begin{enumerate}
\item \underline{\emph{Remove all $2$--bracelets.}} As Figure \ref{fig:bracelets} illustrates, the two boundary components of a $2$--bracelet $B_2$ are isotopic, 
and moreover $B_2$ is homeomorphic to the pair $(\mathbb{S}^2\times I, \{\text{4 points}\}\times I)$.
Thus whenever a $2$--bracelet sits between two other bracelets, those other bracelets can be glued directly to one another, with the gluing map adjusted accordingly.

\smallskip
\item \underline{\emph{Remove needless $1$--bracelets.}} Suppose that a $1$--bracelet $B_1$ is glued to a $d$--bracelet $B_d$ (with $d>1$), and that their preferred slopes at the gluing Conway sphere are Farey neighbors. Then the two arcs of $K_1$ can be isotoped to lie on the Conway sphere $\bdy B_1$, without intersecting the crossing segments of $B_d$. As a result, the arcs $K_d \cup K_1$ combine to form the band of a $(d\!-\!1)$--bracelet, as in Figure \ref{fig:remove-neighbor}. Thus we may remove $B_1$ and replace $B_d$ by a bracelet $B_{d-1}$, with one fewer boundary component.

\begin{figure}[h]
\begin{center}
\psfrag{ar}{$\Rightarrow$}
\includegraphics{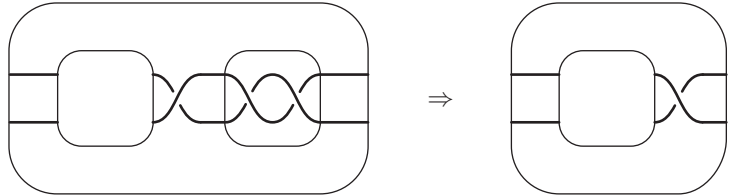}
\caption{Removing a needless $1$--bracelet.}
\label{fig:remove-neighbor}
\end{center}
\end{figure}

\smallskip
\item \underline{\emph{Undo connected sums.}} Suppose that a $1$--bracelet $B_1$ is glued to a $d$--bracelet $B_d$ (with $d>1$), and that their preferred slopes are equal. Then there are several different 2--spheres that pass through the trivial tangle of $B_1$ and intersect $K$ in a pair of points connected by a crossing segment. In this situation, we cut $K$ along the crossing segments of $B_d$, decomposing it as a (possibly trivial) connected sum of $d-1$ other links, as in Figure \ref{fig:undo-connectsum}. On the level of bracelets, each piece of $K$ that was glued to a Conway sphere of $B_d$ will instead be glued to its own $1$--bracelet, whose slope on the gluing sphere is given by the crossing segments of $B_d$.

\begin{figure}[h]
\begin{center}
\psfrag{ar}{$\Rightarrow$}
\psfrag{A}{$A$}
\psfrag{B}{$B$}
\psfrag{C}{$C$}
\includegraphics{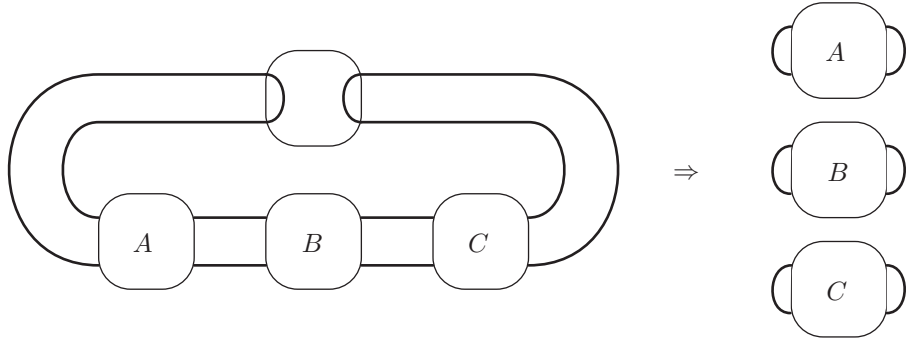}
\caption{Special $1$--bracelets decompose a link as a connected sum.}
\label{fig:undo-connectsum}
\end{center}
\end{figure}

After this cutting operation, we will work separately with each of the $d-1$ new links. That is: we will apply the reduction algorithm to simplify the bracelet presentation of each of these links. The algorithm may reveal that one or more of the new links is actually the unknot (see Theorem \ref{thm:recognition}), and thus that we have undone a trivial connected sum. In this case, we may simply throw away the trivial pieces, having still gained the benefit of a simpler bracelet presentation of $K$.

\smallskip
\item \underline{\emph{Repeat steps $(1)$--$(3)$, as necessary.}} Note that removing a needless $1$--bracelet can create a new $2$--bracelet (as in Figure \ref{fig:remove-neighbor}), and removing a $2$--bracelet can change the gluing map of a $1$--bracelet. However, since each of the above steps reduces the total number of Conway spheres in the construction of $K$, eventually we reach a point where none of these reductions is possible.

\end{enumerate}

\begin{define}\label{def:aug-bracelet}
Let $A \subset \mathbb{S}^3$ be an unknotted band, and let $T = \mathbb{S}^3 \setminus A$ be the open solid torus equal to the complement of $A$. Let $L_n$ be a link consisting of $n$ parallel, unlinked copies of the core of $T$.

Recall from Definition \ref{def:bracelet} that a $d$--bracelet $B_d$ is the pair $(M_d, K_d)$, where $M_d$ is the complement of a regular neighborhood of $d$ disjoint crossing segments of $A$, and $K_d = M_d \cap \bdy A$. We define an \emph{$n$--augmented $d$--bracelet} to be $B_{d,n} = (M_d, K_d \cup L_n)$. Thus a traditional $d$--bracelet $B_d$ corresponds to taking $n=0$. 
When $n\geq 1$, for all positive $d$ (including $d=1$)
the \emph{preferred slope} of $B_{d,n}$ at a boundary (Conway) sphere is the slope of the crossing segments of $A$ at this sphere.
\end{define}

Augmented bracelets naturally arise from certain configurations of $d$--bracelets. We continue our algorithm as follows.

\begin{enumerate}
\setcounter{enumi}{4}

 \item \underline{\emph{Create augmented bracelets.}} Let $B_3$ be a $3$--bracelet glued to $1$--bracelets $B_1$ and $B'_1$. Suppose that there is a marking of the boundary spheres of $B_3$ such that the preferred slope of $B_3$ is $\infty$ and the preferred slopes of $B_1$ and $B'_1$ are in $\mathbb{Z}+1/2$. 
(For example, in the two trivial tangles of Figure \ref{fig:create-augmentation}, these slopes are $-1/2$ and $1/2$. An intrinsic criterion is: the slopes of $B_3$ and $B_1$ (resp. $B'_1$) are not Farey neighbors, but they share exactly two common Farey neighbors.) In this situation, we will replace $B_3 \cup B_1 \cup B'_1$ by a once-augmented 1--bracelet $B_{1,1}$, as in Figure \ref{fig:create-augmentation}. Note that the closed loop in $B_{1,1}$ can be isotoped to lie on the boundary sphere $S$. Up to isotopy, there is exactly one arc pair on $S$ that is disjoint from this loop; its slope is the preferred slope of $B_{1,1}$.

\begin{figure}[h]
\begin{center}
\psfrag{ar}{$\Rightarrow$}
\psfrag{g1}{$\gamma_1$}
\psfrag{g2}{$\gamma_2$}
\includegraphics{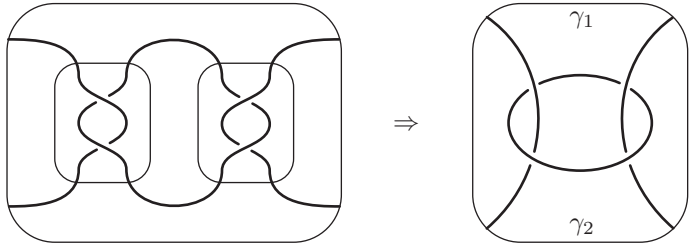}
\caption{Creation of an augmented $1$--bracelet. The arc pair $\gamma_1 \cup \gamma_2$ defines the preferred slope of $B_{1,1}$. Note: if we change the slope in one of the trivial tangles by an integer (e.g. by inverting the two crossings in the left trivial tangle, which will turn its slope $-1/2$ into $+1/2$: a change of $+1$), then the $3$-dimensional picture is the same up to a homeomorphism (a certain number of half--twists in the ``main'' band of the bracelet $B_3$).}
\label{fig:create-augmentation}
\end{center}
\end{figure}

\underline{\emph{Remark:}} 
If $B_3$ is glued to three different $1$--bracelets, each with slope $\pm 1/2$ 
(so the link contains no other bracelets than these four),
we break the symmetry by choosing two of the $1$--bracelets for augmentation. \end{enumerate}

\begin{define}\label{def:large-bracelet}
A (possibly augmented) bracelet $B_{d,n}$ is \emph{large} if $d \geq 3$ or $n \geq 1$. 
\end{define}

\begin{figure}[h]
\begin{center}
\psfrag{ar}{$\Rightarrow$}
\psfrag{p}{$+$}
\includegraphics{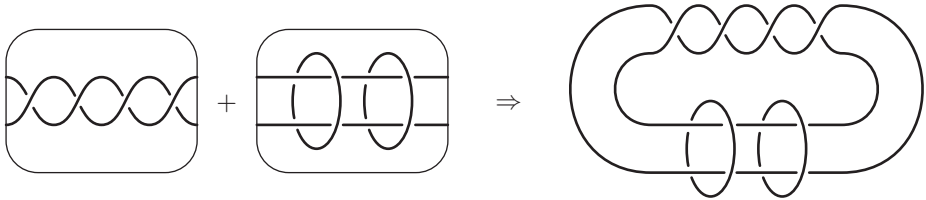}
\vspace{-3ex}
\caption{Creating an augmented $0$--bracelet.}
\label{fig:augment-0bracelet}
\end{center}
\end{figure}

\smallskip
\begin{enumerate}
\setcounter{enumi}{5}
\item \underline{\emph{Combine large bracelets when possible.}} 
Suppose that large bracelets $B_{d,n}$ and $B_{d',n'}$ are glued together along a Conway sphere, with their preferred slopes equal. Then we will combine them into a single $(d+d'-2)$--bracelet, augmented $(n+n')$ times. Note that at the beginning of this step, the only augmented bracelets are of the form $B_{1,1}$, created in step $(5)$. However, under certain gluing maps, several bracelets of this form may combine with other large bracelets to form $n$--augmented $d$--bracelets, with $d$ and $n$ arbitrarily large. 

\smallskip
\item \underline{\emph{Form $0$--bracelets and augmented $0$--bracelets.}} Consider a $1$--bracelet $B_1$, with preferred slope $s$. For any arc pair  $\gamma \subset \bdy B_1$ whose slope is a Farey neighbor of $s$, we can construct a rectangular strip in $B_1$ with boundary $K_1 \cup \gamma$. Therefore, when bracelets $B_1$ and $B'_1$ are glued together and their preferred slopes share a common neighbor in $\ff$, we can glue these two rectangular strips to form an annulus or M\"obius band whose boundary is $K_1 \cup K'_1$. In this situation, we replace $B_1 \cup B'_1$ by a single $0$--bracelet.

In a similar fashion, an augmented $1$--bracelet $B_{1,n}$ contains a rectangular strip whose intersection with the boundary sphere defines the preferred slope of $B_{1,n}$. Therefore, when $B_{1,n}$ is glued to a $1$--bracelet $B_1$ and their preferred slopes are Farey neighbors, we once again have an annulus or M\"obius band. In this situation
(similar to step $(2)$),
we replace $B_1 \cup B_{1,n}$ by a single augmented bracelet $B_{0,n}$, as in Figure \ref{fig:augment-0bracelet}.
\end{enumerate}

%\smallskip

\begin{remark} No further instances of steps $(1)$--$(3)$ occur after the creation of augmented bracelets in step $(5)$. This is clear for step $(1)$: no (unaugmented) $2$--bracelets appear, not even in step $(6)$ because the bracelets that merge in $(6)$ are already large.
Next, observe that the preferred slopes on Conway spheres of large bracelets are never changed after step $(5)$ (not even when Conway spheres are cancelled in $(6)$). An easy discussion then implies that steps $(2)$--$(3)$, or their analogues for augmented bracelets, never occur.
\end{remark}

The following result summarizes the output of the simplification algorithm.

\begin{prop}\label{prop:algorithm-output}
For every generalized arborescent link given as input, the algorithm above produces several ``output links'', of which the input was a connected sum. Let $K$ be such an output link.
Then $K$ is expressed as a gluing of (possibly augmented) bracelets, in which all $2$--bracelets are augmented.

Furthermore, suppose that bracelets $B$ and $B'$ are glued along a Conway sphere. Any path through the $1$--skeleton of the Farey complex connecting the preferred slope of $B$ to the preferred slope of $B'$ must contain at least the following number of edges:

\begin{center}
\begin{tabular}{|r|c|c|}
\hline
{}         & $B'_{1,0}$ & $B'_{d,n}$ large \\ \hline
$B_{1,0}$        & $3$ & $2$ \\ \hline
$B_{d,n}$ large  & $2$ & $1$ \\ \hline
\end{tabular}\end{center}
\end{prop}

\begin{proof}
Observe that the reduction algorithm only changes the topological type of $K$ in step $(3)$, where it cuts $K$ into (possibly trivial) connected summands. Thus the output links do in fact sum to $K$.

Top--left entry of the table: if the preferred slopes of bracelets $B_{1,0}$ and $B'_{1,0}$ are at distance $2$ (or less) in the Farey graph, they share a Farey neighbor and thus step $(7)$ has reduced them to a single $0$--bracelet. Similarly, the reduction of step $(6)$ accounts for the bottom--right entry, and step $(2)$ for the non--diagonal entries.
\end{proof}

\subsection{Analyzing the output}

\begin{figure}[h]
\begin{center}
\psfrag{ar}{$\Rightarrow$}
\includegraphics{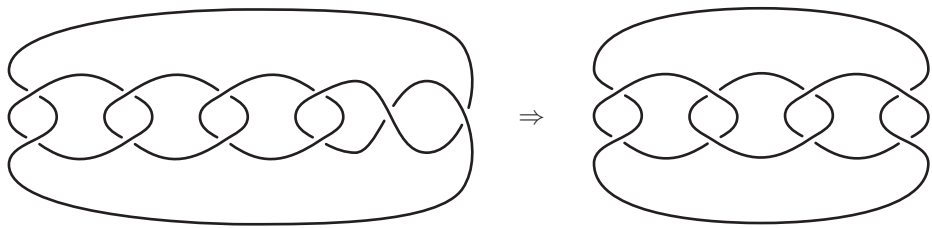}
\caption{The pretzel link $P(2, \minus 2, 2, \minus 2)$ falls in exceptional family II.}
\label{fig:chain-link}
\end{center}
\end{figure}

We are now ready to recognize the non-hyperbolic arborescent links listed in Theorem \ref{thm:main}. After the simplification algorithm, they can appear in any one of four ways:

\begin{enumerate}

\item \underline{\emph{Bracelets augmented more than once.}} A bracelet $B_{d,n}$, where $n \geq 2$, will contain two isotopic link components, as in Figure  \ref{fig:augment-0bracelet}. Each of these parallel components bounds a disk that is punctured twice by the strands of $K_d$. Thus any link containing such a bracelet  falls in exceptional family II.

\smallskip
\item \underline{\emph{$0$--bracelets.}} By Definition \ref{def:bracelet}, the link contained in a $0$--bracelet is the boundary of an unknotted band. These links fall in exceptional family I.

\smallskip
\item \underline{\emph{Once-augmented $0$--bracelets.}} Let $K$ be the link contained in an augmented $0$--bracelet with $r$ half-twists. By reflecting $K$ if necessary, we may assume that $r \geq 0$. Now, we consider three cases:

\begin{itemize}
\item[(a)] \underline{$r=0$}. Then $K$ is the link depicted in Figure \ref{fig:recognition}(a). We note that $K$ is composite, and thus not arborescent by Definition \ref{def:arborescent}.

\quad In fact, we claim that
this case (a) is void because
the reduction algorithm will have cut this link into its prime components
(two copies of the Hopf link).
The augmented $0$--bracelet was necessarily created in step $(7)$ from a $1$--bracelet $B_1$ and an augmented bracelet $B_{1,1}$, which in turn was necessarily created from a $3$--bracelet $B_3$ in step $(5)$. However, $B_3$ must have been glued to $B_1$ with their preferred slopes equal, so in step $(3)$ the algorithm will have recognized $K$ as a connected sum.

\item[(b)] \underline{$r=1$}. Then, as Figure \ref{fig:recognition}(b) shows, $K$ is the boundary of an unknotted band with 4 half-twists, which falls in exceptional family I.

\item[(c)] \underline{$r \geq 2$}. Then, as Figure \ref{fig:recognition}(c) shows, $K$ is the pretzel link $P(r,2,2,\minus 1)$. Because $r \geq 2$ and $\half + \half + \frac{1}{r} > 1$, $K$ falls in exceptional family III.
\end{itemize}

\begin{figure}[h]
\begin{center}
\psfrag{r}{$r$}
\psfrag{a}{(a)}
\psfrag{b}{(b)}
\psfrag{c}{(c)}
\psfrag{ar}{$\Rightarrow$}
\includegraphics{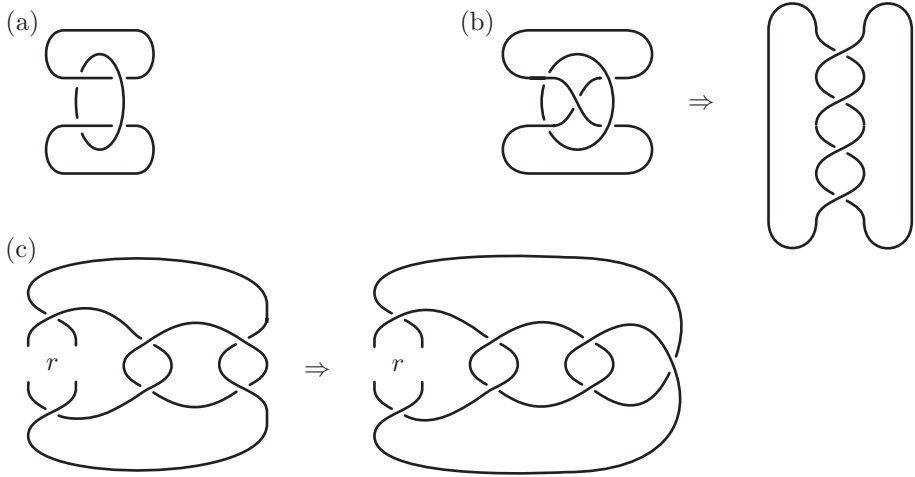}
\caption{Augmented $0$--bracelets form exceptional links in three different ways.}
\label{fig:recognition}
\end{center}
\end{figure}

\smallskip
\item \underline{\emph{Exceptional Montesinos links.}} Recall, from the introduction, that a Montesinos link can be constructed by gluing a bracelet $B_d$ to $d$ different $1$--bracelets. Consider such a link $K$, with $d \geq 3$. By Proposition \ref{prop:algorithm-output}, 
in any Montesinos output link, 
the preferred slope of $B_d$ is not a Farey neighbor of the preferred slope of any of the $1$--bracelets. Thus there is a marking of each Conway sphere, in which the preferred slope of $B_d$ is $\infty$ and the preferred slope of the $1$--bracelet glued to that Conway sphere is not in $\ZZ$.

Once these markings are chosen, there is a unique unknotted band consisting of the arcs of $K_d$ and arcs of slope $0$ along the Conway spheres. We define the \emph{number of half-twists} in the band of $B_d$ to be the number of half-twists in this band. If we modify the marking on some sphere by $k/2$ Dehn twists about the preferred slope of $B_d$, the slope of the $1$--bracelet glued to that sphere goes up by $k$, while the number of half-twists in the band goes down by $k$. Thus, by employing Dehn twists of this sort, we can choose markings in which the preferred slope of $B_d$ is still $\infty$ and the preferred slope of each $1$--bracelet is in the interval $(0,1)$.

We perform two final normalizations. The edge pairs of slopes $0,1,\infty$ decompose each Conway sphere into four triangles. The orientation of the bracelet $B_d$ induces an orientation for each of the $d$ Conway spheres of $\partial B_d$. We stipulate that the edges of any triangle of any Conway sphere $S$ have slope $0,1,\infty$, \emph{in that clockwise order around the triangle}, for the induced orientation of $S$.
(\footnote{We recommend that the reader disregard all orientation issues at a first reading.
The choices in this paragraph are unimportant, as long as they are consistent from one Conway sphere
to the next.})
This can be assumed, up to changing the marking of some of the Conway spheres by an orientation--reversing homeomorphism, exchanging the slopes $0,1$ and fixing $\infty$. Finally, up to reflecting the link $K$, we can ensure that at least one of the $1$--bracelet slopes falls in the interval $(0, 1/2]$. Now, we can recognize the exceptional links:

\begin{itemize}
\item[(a)] If $d=4$, and the slope of every $1$--bracelet is $1/2$, and there are $2$ half-twists in the band of $B_d$, then $K$ is the pretzel link $P(2, \minus 2, 2, \minus 2)$, as in Figure \ref{fig:chain-link}. Thus $K$ falls in exceptional family II.

\item[(b)] If $d=3$, and the slope of every $1$--bracelet is of the form $1/n$, and there is $1$ half-twist in the band of $B_d$, then $K$ is the pretzel link $P(p,q,r,\minus 1)$. If $\frac{1}{p} + \frac{1}{q} + \frac{1}{r} \geq 1$, $K$ falls in exceptional family III.
\end{itemize}
\end{enumerate}

\begin{define}\label{def:candidate}
An arborescent link $K$ is called a \emph{candidate link} if it can be constructed from (possibly augmented) bracelets in the following fashion:
\begin{itemize}
\item $K$ is not a $0$--bracelet or augmented $0$--bracelet.
\item All bracelets are either unaugmented or augmented once.
\item All $2$--bracelets are augmented once.
\item The gluing maps of bracelets along Conway spheres satisfy the minimum--distance table of Proposition \ref{prop:algorithm-output}.
\item If $K$ has exactly one large bracelet $B_{d,0}$, the normalization process of (4a-b) above does not reveal $K$ as an exceptional Montesinos link.
\end{itemize}
\end{define}

So far, we have proved that 
every link $K$ contained in the output of the algorithm (and thus, in particular, every arborescent link)
either falls into one of the three exceptional families, or else is a candidate link. To complete the proof of Theorem \ref{thm:main}, it remains to show the following.

\begin{theorem} \label{thm:sub-main}
Let $K$ be a candidate link. Then $\mathbb{S}^3 \setminus K$ admits a decomposition into angled blocks. Thus, by Theorem \ref{thm:block-hyperbolic}, every candidate link is hyperbolic.
\end{theorem}

The proof of Theorem \ref{thm:sub-main} occupies Sections \ref{sec:complement} and \ref{sec:angles}. In the meantime, we record the following consequence of Theorem \ref{thm:sub-main}, which shows that 
our reduction algorithm does more than merely sort generalized arborescent links into hyperbolic and non-hyperbolic bins.

\begin{theorem}\label{thm:recognition}
Let $K$ be a generalized arborescent link, given in any presentation as a union of unaugmented bracelets. If we apply the simplification algorithm to $K$, it will output several links $K_1, \ldots, K_n$, with the following properties:
\begin{enumerate}
\item $K = K_1 \# \ldots \# K_n$, a connected sum of the output links,
\item every $K_i$ is prime, and
\item if some $K_i$ is the unknot, it appears as a $0$--bracelet with $\pm 1$ half-twist.
\end{enumerate}
Therefore, this algorithm recognizes the unknot and factors a link into its prime summands.
\end{theorem}

\begin{proof}
Statement $(1)$ was part of Proposition \ref{prop:algorithm-output}.
To prove $(2)$, let $K_i$ be a link produced in the output of the algorithm. By Proposition \ref{prop:algorithm-output} and the discussion that follows it, $K_i$ is either a candidate link or a known exception.
By Theorem \ref{thm:sub-main} (which we have yet to prove), the candidate links are hyperbolic, and thus prime. 
As we mentioned in the introduction, the exceptional links in families I and III are also prime. When
$K_i$ is in exceptional family II, we can identify every group of isotopic components to a single circle, producing a new link $K'_i$ that is either hyperbolic or in families I or III. By the cases already discussed, $K'_i$ is prime, and thus $K_i$ is prime also.

To prove $(3)$, suppose that $K_i$ is the unknot. Again, it is easy to check that apart from the $0$--bracelet with $\pm 1$ half-twist, all the exceptional knots discussed above are non-trivial. By Theorem \ref{thm:sub-main}, the candidate links are hyperbolic --- hence also non--trivial. Thus the unknot $K_i$ can in fact be recognized as a $0$--bracelet with $\pm 1$ half-twist. 
\end{proof}

\begin{remark} \label{rem:chain-link}
To solve the general link isotopy problem for arborescent links, Bonahon and Siebenmann \cite{bonsieb-monograph} developed a special type of calculus on weighted trees representing links. Our algorithm for simplifying bracelet presentations is directly inspired by their algorithm for simplifying trees (the latter only keeps more careful track of the order of Conway spheres along each band). The main result of \cite{bonsieb-monograph}, which uses the machinery of double branched covers and equivariant JSJ decompositions, can be paraphrased as follows: for every generalized
arborescent link (with a handful of exceptions), any two bracelet presentations produce the same output under the algorithm above, up to the number of trivial components in the connected sum. The few exceptions include for instance the pretzel link $P(2, \minus 2, 2, \minus 2)$ depicted in Figure \ref{fig:chain-link}, which can appear either in pretzel form (as a $4$--bracelet glued to four $1$--bracelets), or as a twice--augmented $0$--bracelet (in two different ways). Thus, by the results of \cite{bonsieb-monograph}, the algorithm can be said to essentially classify all arborescent links.
\end{remark}

\section{Block decomposition of the link complement}\label{sec:complement}

In this section, we consider a candidate link $K$ (see Definition \ref{def:candidate}) and 
devise a block decomposition of $\mathbb{S}^3\smallsetminus K$ 
(up to applying the algorithm of Section \ref{sec:algorithm}, we may and will assume that the bracelet presentation of $K$ satisfies the conditions of Definition \ref{def:candidate}). In what follows,
we will use one ``large'' block for each large bracelet $B_{d,n}$ ($d\geq 3$ or $n=1$). The Farey combinatorics involved in gluing these blocks to one another and to $1$--bracelets will be encoded in a certain number of ideal tetrahedra in the decomposition. 

Consider a large $d$--bracelet $B_{d,n}$ (where $n=0$ or $1$), and recall that its underlying space $M_d$ is $\mathbb{S}^3$ minus $d$ open balls. Denote by $K_{d,n}$ the tangle (union of arcs) contained in $B_{d,n}$. Between any pair of consecutive Conway spheres, $M_d$ contains a rectangular strip consisting of crossing segments; we call this strip a \emph{crossing rectangle}. The boundary of each crossing rectangle consists of two arcs of $K_d$ and two arcs on Conway spheres, which define the preferred slopes of these two Conway spheres.

The large block $\mathcal{P}_{d,n}$ corresponding to the large bracelet $B_{d,n}$ will be constructed in Section \ref{sec:braceletblocks}: $\mathcal{P}_{d,n}$ will be a solid torus whose boundary is decomposed into ideal polygons, and some pairs of edges (or even faces) of $\mathcal{P}_{d,n}$ will be identified in the block decomposition of $\mathbb{S}^3\smallsetminus K$ (in addition, the core curve of the solid torus $\mathcal{P}_{d,n}$ will be removed if $n=1$). More precisely, each crossing rectangle in $\mathbb{S}^3\smallsetminus K$ will eventually be collapsed to a crossing edge, which will be realized as an ideal edge of $\mathcal{P}_{d,n}$ (in fact, as a pair of identified edges of $\mathcal{P}_{d,n}$ --- see Figure \ref{cylindre} for a preview). Because of this collapsing operation, we need to be very careful that we do not change the nature of the manifold $\mathbb{S}^3\smallsetminus K$ up to homeomorphism. This will be proved in due course (Section \ref{contraction}) before actually constructing $\mathcal{P}_{d,n}$ (Section \ref{sec:braceletblocks}). Meanwhile, we simply insist that the arc pairs of preferred slope (defined e.g. by crossing rectangles) in a Conway sphere of a large bracelet will be realized as ideal edges of the corresponding large blocks. This fact motivates the whole construction of Section \ref{sec:fareytetrahedra}, where we describe all the small blocks (ideal tetrahedra).

\subsection{Gluings through tetrahedra} \label{sec:fareytetrahedra}

In this section, we use the Farey complex $\mathcal{F}$ in order to encode the gluing homeomorphisms between the boundary spheres of two bracelets
into sequences (layers) of ideal tetrahedra, glued at the interface between large bracelets. This follows \cite[Section 5]{gf-bundle}. We also realize trivial tangles as gluings of ideal tetrahedra, following \cite[Appendix]{gf-bundle}. 

\subsubsection{Gluing of two large bracelets}\label{sec:twolargebracelets}

Consider two large bracelets $B_{d,n}$ and $B_{d',n'}$ glued to one another along a Conway sphere.
If $K_{d,n}$ denotes the union of arcs
contained in $B_{d,n}$, define $\mathcal{C}:=M_d\smallsetminus K_{d,n}$ and $\mathcal{C}':=M_d' \smallsetminus K_{d',n'}$. For the later purpose of realizing (a retract of) $\mathcal{C}$ as a polyhedral block, we assume that all the arc pairs of preferred slope are marked on the Conway spheres of $\partial \mathcal{C}$ (and similarly for $\mathcal{C}'$). We will now enhance these graphs on $\partial \mathcal{C}, \partial \mathcal{C}'$ to filling graphs (Definition \ref{def:dual-graph}), then glue $\mathcal{C}$ to $\mathcal{C}'$ \emph{via} a union of ideal tetrahedra attached to the Conway spheres of $\mathcal{C},\mathcal{C}'$ (thus realizing the gluing homeomorphism between $B_{d,n}$ and $B_{d',n'}$).

Let $s$ and $s'$ be the preferred slopes of $\mathcal{C}$ and $\mathcal{C}'$, respectively (see Definition \ref{def:preferred-slope}). By Proposition \ref{prop:algorithm-output}, we know that $s \neq s'$. The obstruction to gluing $\mathcal{C}$ to $\mathcal{C}'$ directly is that arcs of slopes $s$ and $s'$ may have high intersection number. We regard $s,s'\in\mathbb{P}^1\mathbb{Q}$ as vertices in the Farey diagram.

Consider the simplest case, where $s$ and $s'$ are Farey neighbors. Then the four arcs of slope $s$ and $s'$ define a subdivision of a Conway sphere of $\bdy \mathcal{C}$ into two ideal squares. Similarly the arcs of slope $s'$ and $s$ subdivide a Conway sphere of $\partial \mathcal{C}'$ into two ideal squares. In this special case, we may glue  $\bdy \mathcal{C}$ directly to $\partial \mathcal{C}'$ along these two squares.

\begin{figure}[h]
\begin{center}
\includegraphics{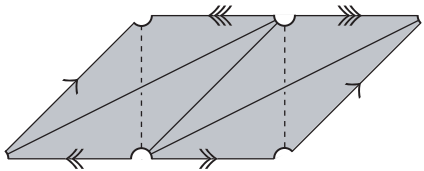}
\caption{A layer of two tetrahedra, caught between two topological pleated surfaces $S_i,S_{i+1}$. Edges with identical arrows are identified.}
\label{fig:tetrahedron-layer}
\end{center}
\end{figure}

If $s,s'$ are not Farey neighbors, we need to consider the sequence of Farey triangles $(T_0,\dots, T_m)$ crossed by the geodesic line from $s$ to $s'$ (where $m\geq 1$, and $s$ [resp. $s'$] is a vertex of $T_0$ [resp. $T_m$]). For each $0\leq i \leq m$, the vertices of $T_i$ define $3$ slopes, and the corresponding arc pairs provide an ideal triangulation $\tau_i$ of the Conway sphere $S$. Moreover, if $x,y$ are the ends of the Farey edge $T_i\cap T_{i+1}$, the two arc pairs whose slopes are $x$ and $y$ define a subdivision $\sigma$ of the $4$--punctured sphere $S$ into two ideal squares: both triangulations $\tau_i$ and $\tau_{i+1}$ are refinements of $\sigma$. In fact, $\tau_{i+1}$ is obtained from $\tau_i$ by a pair of \emph{diagonal exchanges}: remove two opposite edges of $\tau_i$ (thus liberating two square cells, the cells of $\sigma$, of which the removed edges were diagonals); then insert back the other diagonals. Each of these two diagonal exchanges defines (up to isotopy) a topological \emph{ideal tetrahedron} in the space $S\times I$: more precisely, the union of these two ideal tetrahedra is bounded by two topological \emph{pleated surfaces} $S_i,S_{i+1}$ isotopic to $S\times \{*\}$ in $S\times I$, which are pleated along the edges of $\tau_i$ and $\tau_{i+1}$, and intersect each other precisely along the edges of $\sigma$ (see Figure \ref{fig:tetrahedron-layer}). Doing the same construction for all $0\leq i <m$, we thus obtain an ideal triangulation of a strong deformation retract of $S\times I$, whose bottom (resp. top) is pleated along an ideal triangulation of $S$ containing the arc pair of slope $s$ (resp. $s'$).

Finally, the remaining two pairs of edges in the triangulation $\tau_0$ (in addition to the pair of slope $s$) define a subdivision of the boundary of the space $\mathcal{C}$. The same occurs for $\mathcal{C}'$. We have completed our aim of gluing $\mathcal{C}$ to $\mathcal{C}'$, with boundaries suitably triangulated, using a sequence of (pairs of) ideal tetrahedra as an interface. Note that the choice of ``suitable'' boundary triangulations of $\mathcal{C},\mathcal{C}'$ is forced by the gluing homeomorphism itself.

\begin{define} \label{def:product_region}
The family of ideal tetrahedra between $\mathcal{C}$ and $\mathcal{C}'$ is called a \emph{product region}. The same term also refers to the union of that family. Topologically, the product region is a strong deformation retract of $S \cross I$; when $s$ and $s'$ have no common Farey neighbors, the product region is homeomorphic to $S \cross I$.
\end{define}

\subsubsection{Gluing a large bracelet to a trivial tangle} \label{sec:swgluing}

Consider a gluing of a large bracelet $B_{d,n}$ to a trivial tangle $B_{1,0}$ along a Conway sphere. As before, define $\mathcal{C}:=M_d\smallsetminus K_{d,n}$.
and $\mathcal{C}':=M_1\smallsetminus K_{1,0}$. 
We will triangulate the Conway sphere $S$ of $\mathcal{C}$ and attach ideal tetrahedra to $S$ to realize a space homeomorphic to $\mathcal{C} \cup_S \mathcal{C}'$. (Note: we will not need to attach a copy of $\mathcal{C}'$ itself, only ideal tetrahedra.)

Let $s$ and $s'$ be the preferred slopes of $\mathcal{C}$ and $\mathcal{C}'$, respectively. By the minimum--distance table of Proposition \ref{prop:algorithm-output}, we know that $s$ and $s'$ are not equal and are not Farey neighbors. We can thus consider the sequence of Farey triangles $(T_0,\dots,T_m)$ from $s$ to $s'$, where $m\geq 1$. We perform the same construction as in \ref{sec:twolargebracelets} above, using topological pleated surfaces $S_0,\dots, S_{m-1}$ whose triangulations are given by $T_0,\dots, T_{m-1}$ (note the omission of $T_m$). To realize the trivial tangle complement $\mathcal{C}'$, we will now glue the faces of the pleated surface $S_{m-1}$ together in pairs, following Sakuma and Weeks' construction in \cite{sakuma-weeks}.

\begin{figure}%[ht]
\begin{center}
%	\psfrag{gl}{gluing}
%	\psfrag{ar}{$\Longrightarrow$}
%	\psfrag{ad}{$\Downarrow ~~~~\text{isotopy}$}
%	\psfrag{al}{$\Longleftarrow$}
%	\includegraphics{Fig_fold-clasp.eps}
\begin{overpic}{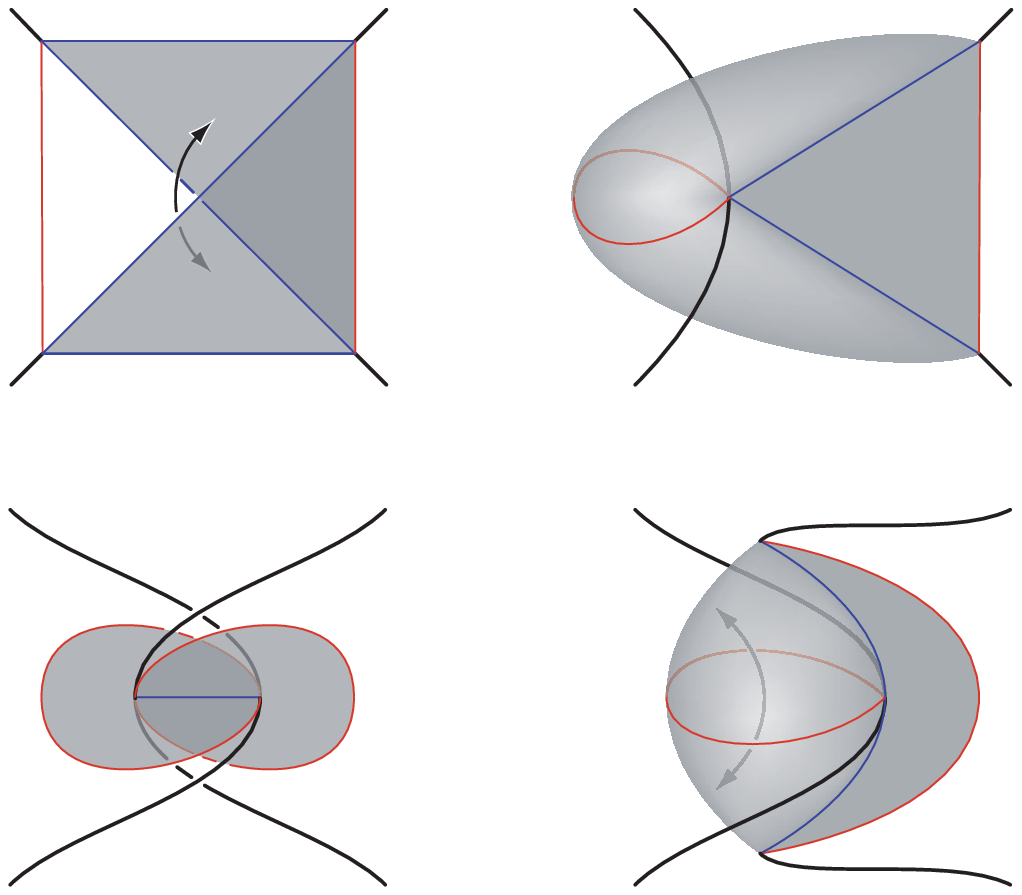}
\put(72, 44){$\Downarrow ~~~~\text{isotopy}$}
\put(44, 70){$\Longrightarrow$}
\put(44, 20){$\Longleftarrow$}
\put(42.5, 66){gluing}
\put(42.5, 16){gluing}
\end{overpic}
\caption{The surface $S_{m-1}$ is glued to itself, realizing a $1$--bracelet (trivial tangle).}
\label{fig:fold-clasp}
\end{center}
\end{figure}

Without loss of generality, we may assume that the vertices of $T_{m-1}$ and $T_m$ are $0,1,\infty$ and $0,1,\frac{1}{2}$ respectively. Each face (ideal triangle) $f$ of $S_{m-1}$ has exactly one edge $e$ of slope $\infty$, shared with another face $f'$. We simply identify $f$ and $f'$ by a homeomorphism respecting $e$. The result is shown in Figure \ref{fig:fold-clasp}: it is straightforward to check that the simple closed curve in $S_{m-1}$ of slope $s'=\frac{1}{2}$ becomes contractible. The gluing thus realizes a $1$--bracelet of slope $s'=\frac{1}{2}$.

\begin{remark} \label{rem:manifoldness}
If $m=1$, note that the Conway sphere of $\mathcal{C}$, made of four ideal triangles, \emph{is} $S_{m-1}$ and has been directly collapsed to two ideal triangles, without gluing any tetrahedra. More generally, for any $m$, all $4$ edges whose slopes are Farey neighbors of $s'$ are collapsed to just one edge (the horizontal edge in the last panel of Figure \ref{fig:fold-clasp}): though none of these four edges can have slope $s$ (because $s,s'$ are not Farey neighbors), some of them may belong to the Conway sphere in $\partial \mathcal{C}$. In spite of these collapsings, for any candidate link $K$, we can realize the space $\mathbb{S}^3\smallsetminus K$ with a well-defined \emph{manifold} structure by gluings of the type above, because the Conway spheres of $\mathcal{C}$ are pairwise disjoint.
\end{remark}

\subsubsection{Gluing two trivial tangles together}

Finally, when two trivial tangles are glued to one another, we obtain a $2$--bridge link $K$. The strands in each bracelet can be isotoped to proper pairs of arcs in the Conway sphere, of slopes $s$ and $s'$. If $s,s'$ are sufficiently far apart in the Farey diagram, we can perform the gluing operation above both near $\mathcal{C}$ and near $\mathcal{C}'$ (both $\mathcal{C}$ and $\mathcal{C}'$ being homeomorphic to $1$--bracelets $M_1 \smallsetminus K_{1,0}$). The resulting decomposition into tetrahedra was constructed by Sakuma and Weeks \cite{sakuma-weeks}, and also described in the Appendix to \cite{gf-bundle}. For completeness, we include

\begin{prop} \label{prop:2bridge-manifoldness}
If $s,s'$ have no common Farey neighbors (i.e. satisfy the mini\-mum--distance table of Proposition \ref{prop:algorithm-output}), the union of the tetrahedra defined by the construction above is a triangulated \emph{manifold} homeomorphic to $\mathbb{S}^3\smallsetminus K$, where $K$ is a $2$-bridge link.
\end{prop}
\begin{proof}
First, the path of Farey triangles $(T_0,\dots,T_m)$ from $s$ to $s'$ satisfies $m\geq 3$: indeed, if $m=2$, then two vertices of $T_1$ are Farey neighbors of $s$, and two are Farey neighbors of $s'$ --- so $s$ and $s'$ have a common neighbor. Therefore the first and last pleated surfaces $S_1$ and $S_{m-1}$ are distinct, and there is at least one layer of tetrahedra. Consider the union of all tetrahedra \emph{before} folding $S_1$ and $S_{m-1}$ onto themselves. Denote by $x,y$ the ends of the Farey edge $T_1\cap T_2$, and thicken the corresponding tetrahedron layer between $S_1$ and $S_2$ by replacing each edge whose slope is $x$ or $y$ with a bigon. The resulting space is homeomorphic to $S_1\times [0,1]$: therefore, after folding $S_1$ and $S_{m-1}$, we do obtain the manifold $\mathbb{S}^3\smallsetminus K$. It remains to collapse the four bigons back to ordinary edges, without turning the space into a non-manifold.

Recall (Remark \ref{rem:manifoldness}) that the folding of $S_{1}$ (resp. $S_{m-1}$) identified all $4$ edges whose slopes are ends of $T_0\cap T_1$ (resp. $T_m\cap T_{m-1}$) to one edge, and caused no other edge identifications. At most one of $x,y$ belongs to the Farey edge $T_0\cap T_1$; at most one of $x,y$ belongs to $T_m\cap T_{m-1}$; and none of $x,y$ belongs to both Farey edges simultaneously ($s,s'$ have no common neighbors). So under the folding of $S_1$ and $S_{m-1}$, the two bigons of slope $x$ may become glued along one edge; the two bigons of slope $y$ may become glued along one edge, and no further identifications occur between points of the $4$ bigons. When two bigons are identified along one edge, consider their union as just one bigon. All (closed) bigons are now disjointly embedded in $\mathbb{S}^3\smallsetminus K$, so we can collapse each of them to an ideal segment without changing the space $\mathbb{S}^3\smallsetminus K$ up to homeomorphism.
\end{proof}

Thus, by the results of \cite{gf-bundle} (Theorem A.1), or more broadly by Menasco's theorem \cite{menasco-alt}, we already know that $2$--bridge links admit angle structures (and are hyperbolic) if and only if they are candidate links.

\subsection{Collapsing} \label{contraction}

At the beginning of Section \ref{sec:complement}, we defined the \emph{crossing rectangles}: a large bracelet $B_{d,n}$ has $d$  crossing rectangles $R\simeq [0,1]\times [0,1]$, such that two opposite sides $\{0,1\}\times [0,1]$ of $R$ define the preferred slopes in two consecutive Conway spheres, and the two other sides $[0,1]\times\{0,1\}$ belong to the tangle (union of arcs) contained in $B_{d,n}$.

In a candidate link $K$ containing at least one large bracelet, we collapse each crossing rectangle $R\simeq [0,1]\times [0,1]$ as above to a segment $\{*\}\times [0,1]$.

\begin{prop}
The space obtained after collapsing the crossing rectangles to segments is still homeomorphic to the manifold $\mathbb{S}^3\smallsetminus K$.
\end{prop}

\begin{proof}
As in Proposition \ref{prop:2bridge-manifoldness}, it is enough to check that the closed crossing rectangles, before collapsing, are disjointly embedded in the union of ideal tetrahedra and uncollapsed large bracelets. First, consider a gluing between two large bracelets: since the two preferred slopes on the gluing Conway sphere are distinct (by the minimum--distance table of Proposition \ref{prop:algorithm-output}), no points of the adjacent crossing rectangles adjacent to this Conway sphere get identified (in the product region corresponding to the Conway sphere, all tetrahedron edges are disjoint). Then, consider a gluing between a large bracelet (with preferred slope $s$) and a trivial tangle (with preferred slope $s'$): in Remark \ref{rem:manifoldness}, we observed that none of the edges which undergo identifications have slope $s$, because $s,s'$ are not Farey neighbors (being at distance $2$ or more in the Farey graph). Therefore, no points of any crossing rectangles are identified. Since the (closed) crossing rectangles are disjointly embedded in $\mathbb{S}^3 \smallsetminus K$, we can collapse each of them to a segment without changing the space $\mathbb{S}^3\smallsetminus K$ up to homeomorphism.
\end{proof}

\subsection{Blocks associated to large bracelets} \label{sec:braceletblocks}
In this section, we construct the blocks associated to large bracelets.
We begin by considering a large \emph{non-augmented} bracelet $B_d$ (where $d\geq 3$) and set out to construct an ideal polyhedron version of the space $\mathcal{C}$, now defined as $M_d\smallsetminus K_d$ with crossing rectangles collapsed to (ideal) segments. We will construct $\mathcal{C}$ as a polyhedral solid torus (or block) $\mathcal{P}$ with some edge identifications.

\begin{figure}[h]
\begin{center}
%	\psfrag{g1}{$\dot{\gamma}_1$}
%	\psfrag{h1}{$\ddot{\gamma}_1$}
%	\psfrag{g2}{$\dot{\gamma}_2$}
%	\psfrag{h2}{$\ddot{\gamma}_2$}
%	\psfrag{a1}{$\gamma_1$}
%	\psfrag{a2}{$\gamma_2$}
%	\includegraphics{Fig_cylindre.eps}
%	\begin{center}
\begin{overpic}{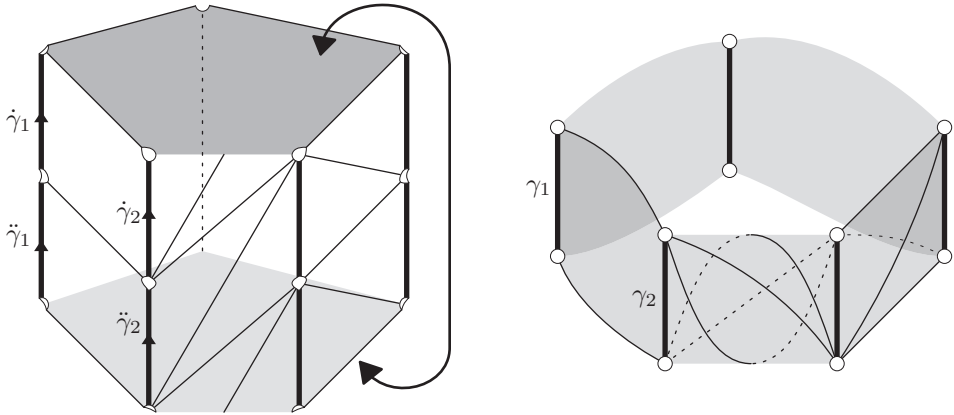}
\put(-3, 32){$\dot{\gamma}_1$}
\put(-3, 19.5){$\ddot{\gamma}_1$}
\put(9, 22){$\dot{\gamma}_2$}
\put(9, 10){$\ddot{\gamma}_2$}
\put(53.5, 25){$\gamma_1$}
\put(65, 13){$\gamma_2$}
\end{overpic}
\end{center}
\vspace{-1ex}
\caption{Block associated to a large (non--augmented) bracelet of degree $5$. The solid torus $\mathcal{P}$ is shown before and after the identifications $\dot{\gamma}_i\sim \ddot{\gamma}_i$. In the right panel, $\mathcal{P}$ is the space \emph{outside} the grey ``ring''.} \label{cylindre}
\end{figure}

Consider a (closed) solid torus $\hat{\pp}$ with a preferred, core--parallel, simple closed curve $\hat{\gamma}$ on $\partial \hat{\pp}$. We endow both $\hat{\mathcal{P}}$ and $\hat{\gamma}$ with orientations that will remain fixed throughout the paper.
Draw $d$ disjoint oriented curves $\hat{\gamma}_1,\dots,\hat{\gamma}_d$ isotopic to $\hat{\gamma}$ in $\partial \hat{\pp}$, so that $\partial \hat{\pp}$ is the union of $d$ annuli $\hat{U}_i$, each bounded by two curves, $\hat{\gamma}_i$ and $\hat{\gamma}_{i+1}$ (indices are taken modulo $d$). For matters of orientation, we assume that the orientation entering $\hat{U}_i$ (from a point of $\hat{\gamma}_i$), followed by the orientation of $\hat{\gamma}_i$, is equal to the orientation on $\partial \hat{\mathcal{P}}$ induced by the orientation of $\hat{\mathcal{P}}$ (Figure \ref{cylindre}, left). Next, remove one point in each $\hat{\gamma}_i$. We still denote by $\hat{\pp}$ the solid torus with these $d$ points removed, and by $\hat{U}_i$ the annuli with two boundary points removed. Each curve $\hat{\gamma}_i$ is now replaced by an ideal (oriented) arc, also denoted $\hat{\gamma}_i$, joining a puncture to itself.

To construct an ideal cellulation of $\partial \hat{\pp}$, we must decompose each annulus $\hat{U}_i$ into contractible ideal polygons. There are two options for doing so:
\begin{enumerate}
\item Choose an ideal arc across $\hat{U}_i$, connecting the puncture of $\hat{\gamma}_i$ to the puncture of $\hat{\gamma}_{i+1}$ (for each $i$, there is a $\mathbb{Z}$--worth of possible choices for such an arc). Then $\hat{U}_i$ is an ideal square cell with a pair of opposite sides identified.
\item Choose two disjoint, non-isotopic ideal arcs across $\hat{U}_i$, decomposing $\hat{U}_i$ into two ideal triangles. There is again a $\mathbb{Z}$--worth of possible choices.
\end{enumerate}
Finally, consider the two-fold cyclic cover $\pp$ of $\hat{\pp}$. Each $\hat{\gamma}_i$ lifts to two arcs $\dot{\gamma}_i$ and $\ddot{\gamma}_i$ in $\partial \pp$ with the same pair of (distinct) punctures as end points. 
We now identify $\dot{\gamma}_i$ with $\ddot{\gamma}_i$, by an orientation--reversing homeomorphism. The resulting arc is called $\gamma_i$, and the quotient of $\pp$ under this identification is homeomorphic to $\mathcal{C}$, as defined at the beginning of this section. 
(The $\gamma_i$ are the collapsed crossing rectangles, and the removed tangle $K_d$ is ``in the ideal vertices''.) Each arc $\gamma_i$ is a crossing segment, not endowed with any orientation. See Figure \ref{cylindre}.

Under the identification, each band $U_i$ (the lift of $\hat{U}_i$ in $\pp$) becomes a $4$--punctured sphere included in the boundary of $\mathcal{C}$. Any decomposition of $\hat{U}_i$, as in the dichotomy above, defines a decomposition of $U_i$ into two ideal squares or four ideal triangles.
The two cases of the dichotomy correspond, as in Section \ref{sec:fareytetrahedra}, to the neighboring bracelet either $(1)$ being glued directly to $\mathcal{C}$, using a very ``simple'' homeomorphism, or $(2)$ needing an interface of ideal tetrahedra or being in fact a $1$--bracelet (trivial tangle).

Recall that in each case of the dichotomy, the involved choices in $\mathbb{Z}$ are forced by the gluing homeomorphisms (see Section \ref{sec:fareytetrahedra}), and essentially reflect the number of half-twists in the band of the bracelet $(M_d, K_d)$ that defined $\mathcal{C}$ (and the block $\mathcal{P}$).

\subsubsection{Augmented polyhedral $d$--bracelets ($d\geq 1$)} \label{sec:augmentedblocks}

Topologically, the construction above is valid for all $d\geq 1$, not just $d\geq 3$; and an augmented $d$--block (associated to an augmented $d$--bracelet) is obtained by drilling out the core of the solid torus $\pp$. If we denote by $\pp_{aug}$ the result of the drill-out, then any decomposition of $\partial \pp$ into ideal triangles and squares (as in the dichotomy above) induces, by 
coning off to the core of $\mathcal{P}_{aug}$, a decomposition of $\pp_{aug}$ into contractible ideal polyhedra (tetrahedra and square--based pyramids): the coning--off is induced by the product structure $\mathbb{T}^2\times (0,1)$ of the interior of $\mathcal{P}_{aug}$.
For the purpose of finding angle structures, it will be convenient to regard $\mathcal{P}_{aug}$ as such a union of ideal polyhedra, rather than an elementary block \emph{per se}.

An essential feature of our polyhedral realizations of $d$--blocks and augmented $d$--blocks is that crossing segments always arise as edges, thus defining a preferred slope on each Conway sphere of the block. If the block comes from an augmented bracelet, this slope is also the one defined by the ``extra'' link component living inside the augmented bracelet.

\section{Angle structures for the link complement}\label{sec:angles}

In this section, we find dihedral angles (satisfying Definition \ref{def:angled-block} and the hypotheses of Theorem \ref{thm:block-hyperbolic})
for the blocks and ideal tetrahedra constructed in Section \ref{sec:complement}. The key property will be an explicit description (Section \ref{sec:concave-hops}) of the space of angle structures associated to a sequence of tetrahedra forming a product region (see Definition \ref{def:product_region}) or a $1$--bracelet (trivial tangle). This description (in large part borrowed from \cite[Section 5 \& Appendix]{gf-bundle}) is sufficiently tractable that we can say exactly when the tetrahedra admit dihedral angles that match a given system of angles for the solid tori (Propositions \ref{prop:innerbranch} and \ref{prop:paqb}). 
As a result, we can show that all candidate links (see Definition \ref{def:candidate}) admit angle structures. We will treat the easier case of non-Montesinos links in Section \ref{sec:non-montesinos} and the trickier case of Montesinos links in Section \ref{sec:montesinos}. Montesinos links are tricky because they include the third family of exceptions to Theorem \ref{thm:main}.

The strategy is as follows: in Section \ref{sec:angledtorus}, we choose some dihedral angles to parametrize the deformation space of solid tori. We show (Propositions \ref{prop:chi-angles}, \ref{prop:degreed}, \ref{prop:degree3}) that these dihedral angles define valid angle structures whenever they are, in some appropriate sense, small enough. On the other hand, the same parameters need to be large enough (Proposition \ref{prop:paqb}) for the trivial tangles to admit angle structures. The conflict that arises can be managed for candidate links, but causes the exceptional Montesinos links (ruled out in Definition \ref{def:candidate}) to have no angle structures.

\subsection{Angle structures for a non-augmented large block}\label{sec:angledtorus}

We consider a solid torus $\pp$ of degree $d$ whose boundary is subdivided into ideal triangles and quadrilaterals, as in Section \ref{sec:braceletblocks}. In this section, we study the space of angle structures for $\pp$.
We restrict attention to those angle structures which are invariant under the natural fixed-point-free involution of $\pp$ (recall that $\pp$ was defined as a $2$-fold covering in Section \ref{sec:braceletblocks}).

Recall the preferred direction of $\partial \pp$, defined by the crossing arcs of the corresponding bracelet. As parameters, we will use the (exterior) dihedral angles at all those edges of $\partial \pp$ which are \emph{not} along the preferred direction. The angles at the edges along the preferred direction can then be recovered from the requirement that the angles around any vertex of the block add to $2\pi$.

For simplicity, assume that $\partial \pp$ is decomposed into ideal triangles only. Recall from Section \ref{sec:braceletblocks} that the preferred direction came with an orientation, which we call ``upwards''. Thus, each of the $d$ bands $U_1,\dots,U_d$ in $\partial \pp$ is traversed by two \emph{ascending} and two \emph{descending} edges
(e.g. in Figure \ref{cylindre} [left], though all four edges across $U_2$ seem to go upwards to the right, we agree to call only the steeper pair ascending, and the other pair descending).
By the normalization of markings that precedes Definition \ref{def:candidate}, the slopes of the descending edges in the corresponding Conway spheres are $0$, while the slopes of the ascending edges are $1$, and the slopes associated to the edges along the preferred direction are $\infty$.

For each $1\leq i \leq d$, we denote by 
a formal variable
$a_i \in [0,\pi)$ the (exterior) dihedral angle at the ascending pair of edges, and $b_i\in [0,\pi)$ the angle at the descending pair. (We can recover the case where $U_i$ is subdivided into two squares by allowing $a_i=0$ or $b_i=0$). If the index $i$ is read modulo $d$, the exterior dihedral angle at an edge (in the preferred direction) of $U_i\cap U_{i+1}$ must be
\begin{equation}\label{eq:ciip}
c_{i,i+1}=\pi-\frac{a_i+a_{i+1}+b_i+b_{i+1}}{2}~,
\end{equation}
so that the angles around each ideal vertex add to $2\pi$. 
To force all angles to be non-negative, we require
\begin{equation}\label{eq:aibi} 0\leq a_i<\pi~\text{ and }~0\leq b_i<\pi~\text{ and }~0<a_i+b_i\leq \pi~.\end{equation} The particular choice of strong and weak inequalities here implies 

\begin{lemma}
Any triangular face of $\mathcal{P}$ has at most one edge with dihedral angle 0.
\end{lemma}

\begin{proof}
The exterior dihedral angles of a triangular face are $a_i, b_i, c_{i,i+1}$ or $a_i, b_i, c_{i-1,i}$.
If $a_i=0$ one easily checks that $b_i, c_{i,i+1},c_{i-1,i}$ are all positive (because $b_i=a_i+b_i>0$ and $a_i+b_i=b_i<\pi$). If $b_i=0$ the argument is the same. If $c_{i,i+1}=0$ then $a_i,a_{i+1},b_i,b_{i+1}$ are all positive (because their sum is $2\pi$).
\end{proof}

Thus, a block $\mathcal{P}$ with only triangular faces and \emph{non-negative} dihedral angles $a_i, b_i$ satisfying (\ref{eq:aibi}) uniquely defines an block with \emph{positive} dihedral angles and, possibly, some quadrilateral faces. By default, we will usually consider that a block $\mathcal{P}$ has only triangular faces and look for angle systems satisfying (\ref{eq:aibi}). 

From Definition \ref{def:angled-block}, recall that an \emph{angle structure} on $\pp$ requires every normal simple closed curve $\gamma$ in $\partial \pp$ that bounds a disk in $\pp$ to have total bending number larger than $2\pi$ (the total bending number is the sum of the exterior dihedral angles at the edges encountered by $\gamma$, counted with multiplicity). Such a curve $\gamma$ can be defined as a non-backtracking closed path in the dual graph. 

\begin{figure}%[h]
\psfrag{a}{$\mathcal{A}$}
\psfrag{b}{$\mathcal{B}$}
\psfrag{f}{$f$}
\psfrag{g}{$g$}
\psfrag{ui}{$\widetilde{U}_i$}
\psfrag{uj}{$\widetilde{U}_j$}
\psfrag{c}{$\widetilde{\gamma}$}
\begin{center}
\includegraphics{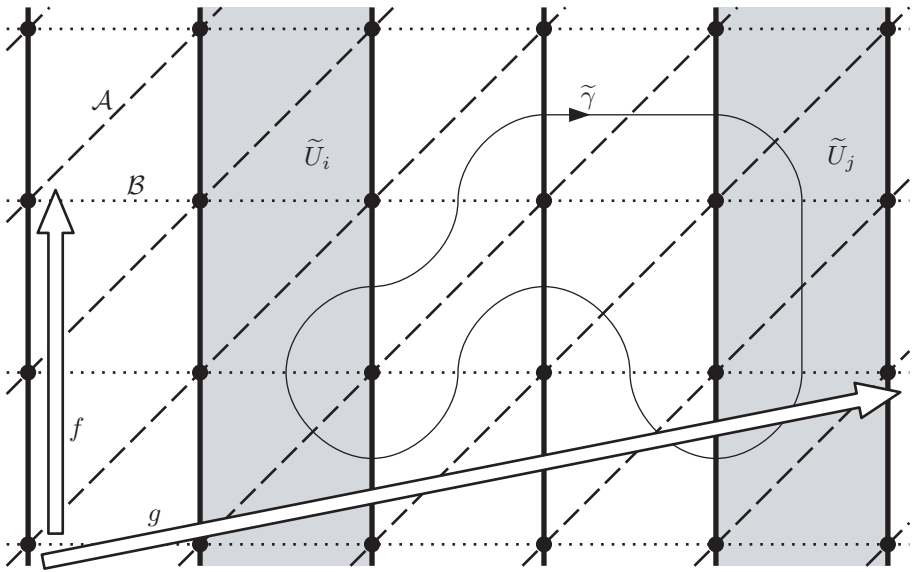}
\end{center}
\vspace{-1ex}
\caption{The cover $\mathbb{R}^2\smallsetminus \mathbb{Z}^2$ of $\partial \mathcal{P}$, with a closed curve $\widetilde{\gamma}$.}
\label{pipeline}
\end{figure}

Identify $\partial \pp$ with $(\mathbb{R}^2\smallsetminus \mathbb{Z}^2)/\langle f,g \rangle$ where $f(x,y)=(x,y+2)$ and $g(x,y)=(x+d,y+k)$, where $d\geq 3$ is the degree of $\pp$ and $k\in \mathbb{Z}$ is an integer such that any $g$-invariant straight line in $\mathbb{R}^2\smallsetminus \mathbb{Z}^2$ projects in $\partial \pp$ to the boundary of a compression disk of $\pp$. 
(In other words, $k$ 
is, up to a constant,
the number of half-twists in the band of the bracelet associated to $\pp$.) The bands $U_i\subset \partial\pp$ lift to a subdivision of $\mathbb{R}^2\smallsetminus \mathbb{Z}^2$ into ``vertical'' bands $(\widetilde{U}_i)_{i\in\mathbb{Z}}$, where $\widetilde{U}_i=(i,i+1)\times \mathbb{R}$ (see Figure \ref{pipeline}). Consider a normal simple closed curve $\gamma$ in $\partial \pp$ bounding a disk in $\pp$, and lift $\gamma$ to a curve $\widetilde{\gamma}$ in $\mathbb{R}^2 \smallsetminus \mathbb{Z}^2$.

Either $\widetilde{\gamma}(1)=\widetilde{\gamma}(0)$ (i.e. $\widetilde{\gamma}$ is a closed curve), or $\widetilde{\gamma}(1)=g^{\pm 1}(\widetilde{\gamma}(0))$ (the exponent cannot be larger than $1$ in absolute value, because $\gamma$ is simple).

\begin{lemma} If $\widetilde{\gamma}(1)=\widetilde{\gamma}(0)$, the total bending number $N_{\gamma}$ of $\gamma$ is more than $2\pi$, unless $\gamma$ is just a loop around a puncture of $\partial \mathcal{P}$.
\end{lemma}

\begin{proof}
If $\widetilde{\gamma}$ visits exactly the vertical bands $\widetilde{U}_i,\widetilde{U}_{i+1},\dots,\widetilde{U}_{j-1},\widetilde{U}_j$, then $\widetilde{\gamma}$ crosses both ascending and descending edges in $U_i$ and in $U_j$ (because $\widetilde{\gamma}$ never crosses the same edge twice consecutively). Counting edges met by $\widetilde{\gamma}$, we thus get

\begin{eqnarray*} N_{\gamma} &\geq& 
(a_i+b_i)+(a_j+b_j)+2\sum_{k=i+1}^{j-1}\min\{a_k,b_k\}+2\sum_{k=i}^{j-1}c_{k,k+1}\\&= & 2\pi+2\sum_{k=i+1}^{j-1}\pi-\max\{a_k,b_k\}~, \end{eqnarray*}
using (\ref{eq:ciip}). This quantity is larger than $2\pi$, unless $j=i+1$ and $\widetilde{\gamma}$ is (homotopic to) the boundary of a regular neighborhood of the union of $r$ consecutive vertical edges (along the preferred direction): in the latter case however, \begin{eqnarray*}N_{\gamma}&=&(r+1)(a_i+b_i)+(r+1)(a_{i+1}+b_{i+1})+2c_{i,i+1}\\ &=&2\pi+r(a_i+b_i+a_{i+1}+b_{i+1})\end{eqnarray*} is larger than $2\pi$, unless $r=0$ and $\gamma$ is just a loop around a puncture of $\partial \pp$. \end{proof}

There remains the case where $\widetilde{\gamma}(1)=g^{\pm 1}(\widetilde{\gamma}(0))$, i.e. $\gamma$ bounds a compression disk of the solid torus $\pp$. Then, the condition $N_{\gamma}>2\pi$ (in terms of the $a_i, b_i$) is in general non-vacuous, and the conjunction of all these conditions (for all normal curves $\gamma$) has no nice closed-form expression in terms of the $a_i, b_i$.
However, the following three Propositions give various sufficient conditions for $N_{\gamma}>2\pi$ to hold, independently of $\gamma$.
\begin{prop} \label{prop:chi-angles}
If the block $\pp$ is of degree $d\geq 3$ and if \begin{equation}\label{eq:chi-angles} (d-2)\pi>\sum_{i=1}^d \max\{a_i,b_i\} \hspace{12pt} \text{\emph{(girth condition),}} \end{equation}
then $N_{\gamma}>2\pi$ for all $\gamma$ bounding a compression disk, i.e. the angles $a_i,b_i,c_{i,i+1}$ define an angle structure for $\pp$.
\end{prop}
\begin{proof}
If $\gamma$ bounds a compression disk in $\pp$, then $\gamma$ crosses all the bands $U_i$. Therefore, $$N_{\gamma}\geq \sum_{i=1}^d\min\{a_i, b_i\}+\sum_{i=1}^d c_{i,i+1}=d\pi-\sum_{i=1}^d \max\{a_i,b_i\}~,$$ using (\ref{eq:ciip}). The conclusion follows.
\end{proof}

Consider a normal curve $\gamma\subset\partial \pp$ bounding a compression disk in $\pp$. 
Let $\mathcal{A}$ (resp. $\mathcal{B}$) be the union of all the ascending arcs (resp. descending arcs) across the annuli $U_i\subset \partial \pp$, each oriented from $U_{i-1}$ to $U_{i+1}$. The torus $\mathbb{T}:=\mathbb{R}^2/\langle f,g\rangle$ naturally contains $\partial \pp$ as a subset, and the closure $\overline{\mathcal{A}}$ of $\mathcal{A}$ (resp. $\overline{\mathcal{B}}$ of $\mathcal{B}$) in $\mathbb{T}$ defines a union of oriented, parallel simple closed curves in $\mathbb{T}$ (see Figure \ref{pipeline}).

\begin{define} \label{def:na-nb}
We denote by $n_{\mathcal{A}}$ (resp. $n_{\mathcal{B}}$) the absolute value of the (algebraic) intersection number of $\gamma$ with $\overline{\mathcal{A}}$ (resp. $\overline{\mathcal{B}}$) in $\mathbb{T}$. This definition clearly does not depend on the choice of compression--disk--bounding curve $\gamma$ (all such $\gamma$ are in the same free homotopy class of $H_1(\mathbb{T},\mathbb{Z})$).
\end{define}

\begin{remark} \label{rem:count-twists}
In the special case of Montesinos links, we defined an integer $n$ (up to sign), called the number of half--twists in the main band, in the normalization that precedes Definition \ref{def:candidate}. Since each descending arc has slope $0$ in the corresponding Conway sphere, the definition of $n$ implies that $|n|=n_{\mathcal{B}}$ for Montesinos links.
\end{remark}

\begin{prop} \label{prop:degreed}
If $n_{\mathcal{A}}\geq 3$ and $(\alpha_i)_{1\leq i \leq d},(\beta_i)_{1\leq i \leq d}$ are positive numbers such that $\alpha_i\geq \beta_i$, then setting $$(a_i, b_i)=(\pi-\varepsilon\alpha_i, \varepsilon\beta_i)$$ defines an angle structure for the block $\pp$ for all sufficiently small $\varepsilon$.
\end{prop}
\begin{proof}
First, the $a_i, b_i$ clearly satisfy Condition (\ref{eq:aibi}) above. Any normal curve $\gamma$ bounding a compression disk in $\pp$ meets at least $n_{\mathcal{A}}\geq 3$ ascending edges, whose pleating angles are all close to $\pi$: thus $N_{\gamma}>2\pi$ for some small enough $\varepsilon$ (independent of $\gamma$). A similar proposition holds when $n_{\mathcal{B}}\geq 3$.
\end{proof}

\begin{prop} \label{prop:degree3}
Suppose the block $\pp$ has degree $d=3$ and $(n_{\mathcal{A}},n_{\mathcal{B}})=(2,1)$. Pick positive numbers $(\alpha_i)_{1\leq i \leq d},(\beta_i)_{1\leq i \leq d}$ such that $\alpha_i\geq \beta_i$. Setting $(a_i, b_i)=(\pi-\varepsilon\alpha_i, \varepsilon\beta_i)$ defines an angle structure on $\pp$ (for small $\varepsilon$) if and only if $$\alpha_i > \beta_{i+1}+\beta_{i-1}$$ for each $i\in\{1,2,3\}$ (taking indices modulo $3$).
\end{prop}
\begin{proof}
Again, Condition (\ref{eq:aibi}) is clearly satisfied.
If $\varepsilon$ is small enough, as in Proposition \ref{prop:degreed}
it is enough to check $N_{\gamma}>2\pi$ for those compression--disk--bounding curves $\gamma$ which cross $\mathcal{A}$ exactly twice. There are only three such curves: each of them crosses two ascending and one descending edge, hence for some $1\leq i \leq 3$ we have
$$N_{\gamma}=a_{i-1}+b_i+a_{i+1}+\sum_{i=1}^3 c_{i,i+1}=2\pi+\alpha_i-(\beta_{i+1}+\beta_{i-1}).$$ The conclusion follows. A similar proposition holds when $(n_{\mathcal{A}},n_{\mathcal{B}})=(1,2)$.
\end{proof}

\subsubsection{Angle structures on augmented blocks}

Finally, we note that an \emph{augmented} block with prescribed non-negative dihedral angles $a_i, b_i$ (where (\ref{eq:aibi}) holds) can always be realized as a union of tetrahedra with positive angles: the space obtained by coning off the band $U_i$ to the extra component of the tangle in the block (as in Section \ref{sec:augmentedblocks}) is a union of $4$ isometric tetrahedra of interior dihedral angles 
\begin{equation} \frac{\pi-a_i}{2}~,~\frac{\pi-b_i}{2}~,~\frac{a_i+b_i}{2}~, \label{eq:aug-angles} \end{equation} 
all positive by (\ref{eq:aibi}). The exterior dihedral angles of the augmented block are recovered as $\pi-\left ( \frac{\pi-a_i}{2} + \frac{\pi-a_i}{2} \right )=a_i~$, similarly $\pi- \left ( \frac{\pi-b_i}{2} + \frac{\pi-b_i}{2} \right )=b_i~$, and $\pi- \left ( \frac{a_i+b_i}{2}+\frac{a_{i+1}+b_{i+1}}{2}\right )=c_{i,i+1}$.

\begin{remark}
The augmentation component of each augmented $d$--bracelet bounds $d$ disjoint, homotopically distinct, twice--punctured disks (also known as thrice--punctured spheres), which must be totally geodesic for the hyperbolic metric if one exists.
\end{remark}

In fact, suppose \emph{all} large blocks are augmented (the candidate link $K$ is called \emph{totally augmented}; 
a special case of this is the case where there are no large blocks, and $K$ is a 2-bridge link). For totally augmented links, 
the triangulation constructed above falls into the class studied in \cite[Chapter 2]{fg-thesis}. There, it was shown that the triangulation not only admits positive angle structures, but that one of these structures (the one with largest volume) actually realizes the hyperbolic metric, and is a refinement of the \emph{geometrically canonical decomposition} in the sense of Epstein and Penner \cite{epstein-penner}. In other words, a certain (explicit) coarsening of the triangulation is combinatorially dual to the Ford--Vorono\"{\i} domain of the manifold with respect to horoball neighborhoods of the cusps which are chosen to be pairwise tangent at each thrice--punctured sphere.

\subsection{Angle structures for product regions and trivial tangles}\label{sec:concave-hops}

In this section, we investigate the space of angle structures for the ideal tetrahedra constructed in Section \ref{sec:fareytetrahedra}. Tetrahedra live either at the interface of large blocks $\pp,\pp'$, or near trivial tangles. While the space of angle structures for a tetrahedron is easy to describe (a triple of positive angles summing to $\pi$), the difficulty is to deal with many tetrahedra simultaneously.

We begin by focusing on two large blocks $\pp,\pp'$ separated by a product region (these come from large bracelets, in the sense of Definition \ref{def:large-bracelet}). 
Recall from Section \ref{sec:fareytetrahedra} the pleated $4$--punctured spheres $S_j$ between $\pp$ and $\pp'$: we can endow $S_j$ with a transverse ``upward'' orientation, from $\pp$ to $\pp'$. Suppose that we \emph{have} solved the problem of finding an angle structure, i.e. that the tetrahedra and solid tori are assigned dihedral angles that add up to $2\pi$ around each edge.
Then we can define the \emph{pleating angle} of the surface $S_j$ at any edge $e\subset S_j$: namely, if the sum of all dihedral angles at $e$ of tetrahedra and/or solid tori above (resp. below) $S_j$ is $\pi+\alpha$ (resp. $\pi-\alpha$), we say that $S_j$ has pleating angle $\alpha$ at $e$. 
 
It will turn out that pleating angles of the $S_j$ are very convenient parameters for the space of angle structures: thus, when no angle structure has been defined yet, we will typically look for angle structures realizing a given set of pleating angles of the $S_j$, and express the dihedral angles of the blocks $\pp,\pp'$ and ideal tetrahedra in terms of these pleating angles.

We arbitrarily require that \emph{the pleating angles of $S_j$ at the $3$ edges adjacent to any puncture of $S_j$ add up to $0$} (note that this property would hold in a true hyperbolic metric, where the holonomy of the loop around any puncture is a parabolic element of $\text{Isom}^+\mathbb{H}^3$). This property easily implies that the pleating angles of $S_j$ at opposite edges are equal. Restricting to such a subspace of solutions might (in principle) hamper our goal of finding angle structures; however, it is technically very convenient, for reasons we are about to outline now.

Consider the $4$--punctured sphere $S^{(i)}$ defined by the vertical band $U_i\subset \partial\pp$ (by identifying the edges in $\partial U_i$ to create the crossing arcs, as in Section \ref{sec:braceletblocks}). The pleating angles of $S^{(i)}$, in the above convention, are $a_i, b_i$ and $-a_i-b_i$, the latter being the angle at the crossing arcs. Similarly, the $4$--punctured sphere $S^{(i+1)}$ defined by the band $U_{i+1}$ (i.e. corresponding to the \emph{next} Conway sphere) has pleating angles $a_{i+1},b_{i+1},-a_{i+1}-b_{i+1}$. Let $e$ be the crossing edge $S^{(i)}\cap S^{(i+1)}$: recall that $e$ is obtained by identification of two edges of the solid torus $\mathcal{P}$, both carrying an interior dihedral angle of $\pi-c_{i,i+1}$. If the interior dihedral angles at $e$ above $S^{(i)}$ (resp. $S^{(i+1)}$) for the transverse orientation add up to $\pi-(a_i+b_i)$ (resp. $\pi-(a_{i+1}+b_{i+1})$), the sum of all dihedral angles at $e$ will be $$[\pi-(a_i+b_i)]+[\pi-(a_{i+1}+b_{i+1})]+2[\pi-c_{i,i+1}]=2\pi~.$$ Therefore the linear gluing equation at $e$ will automatically be satisfied.

Recall the Farey vertices $s\neq s'$ from Section \ref{sec:fareytetrahedra} associated to the crossing arcs (or preferred slopes) of $\pp$ and $\pp'$. If $s,s'$ are Farey neighbors, then $\pp$ and $\pp'$ are glued directly to one another along a $4$--punctured sphere $S$: the edge pairs of slopes $s,s'$ subdivide $S$ into two ideal squares, and the bands in $\partial \pp,\partial \pp'$ are traversed by edges exactly as in the first member, (1) vs (2), of the dichotomy of Section \ref{sec:braceletblocks} (up to a degree $2$ covering). 

\begin{prop} \label{prop:empty-branch} At the two parallel edges traversing the band of $\partial\pp$, we put an (exterior) dihedral angle $\varepsilon>0$. We put the same angle $\varepsilon$ at the edges traversing the band of $\pp'$. Then, the pleating angles of $S$ at the edge pairs of slope $s,s'$ are $-\varepsilon,\varepsilon$ respectively. \qed \end{prop}

The previous proposition is obvious. Moreover, observe that we can artificially select a pair of diagonals in the two squares making up $S$ and define the third pleating angle (along these diagonals) to be $0$: then (\ref{eq:aibi}) is satisfied because the exterior dihedral angles $a_i, b_i$ of $\mathcal{P}$ (resp. $\mathcal{P}'$) at the gluing Conway sphere are $0$ and $\varepsilon$, though not necessarily in that order.

\smallskip

We now consider the case where $s,s'$ are not Farey neighbors. The bands of $\partial \pp,\partial\pp'$ are now subdivided into 4 triangles each (as in the second member of the dichotomy of Section \ref{sec:braceletblocks}), defining a pair of ascending and a pair of descending edges in each of the two bands. Fix an arbitrary marking of the Conway sphere along which $\mathcal{P}$ is glued to $\mathcal{P}'$. Denote by $A\in \mathbb{P}^1\mathbb{Q}$ (resp. $B\in \mathbb{P}^1\mathbb{Q}$) the slope of the ascending (resp. descending) edge pair in the band of $\partial \pp$, and denote similarly by $A',B' \in\mathbb{P}^1\mathbb{Q}$ the slopes of the edges in $\partial \pp'$. We make no assumption on the order of $A,A',B,B'$ in $\mathbb{P}^1\mathbb{Q}$, i.e. we favor no convention as to which pair is ascending and which is descending. Denote by $a,b$ the exterior dihedral angles of the block $\pp$ at the ascending and descending edges respectively, and define $a',b'$ in a similar way (relative to $\mathcal{P}'$).

\begin{figure}[h!] \centering 
\psfrag{w}{$-w_i$}
\psfrag{w1}{$w_{i-1}$}
\psfrag{ww}{$w_i\!-\!w_{i-1}$}
\psfrag{e}{$e_i$}
\psfrag{e1}{$e_{i-1}$}
\psfrag{L}{$L$}
\psfrag{R}{$R$}
\psfrag{a}{$u$}
\psfrag{bb}{$-w$}
\psfrag{b}{$w$}
\psfrag{cc}{$-v$}
\psfrag{cb}{$v\!-\!w$}
\psfrag{ba}{$w\!-\!u$}
%
%   FG The following are labels for the vertices of the Farey triangles.
%
\psfrag{p}{}%{$p$}
\psfrag{p1}{}%{$p'$}
\psfrag{p2}{}%{$p''$}
\includegraphics{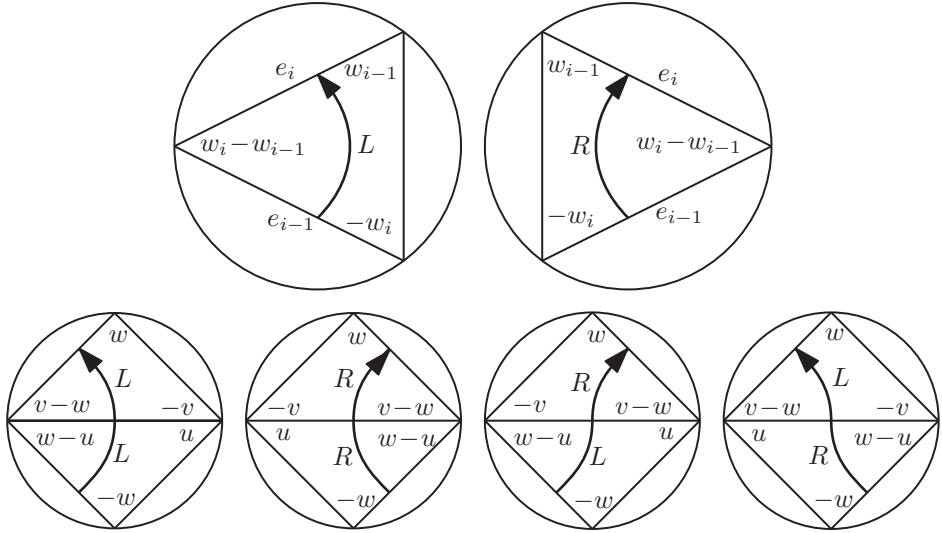}
\caption{Bottom: $e_i$ is the horizontal edge and $(w_{i-1},w_i,w_{i+1})=(u,w,v)$.} \label{fig:rl2}
\end{figure}

\begin{prop} \label{prop:innerbranch}
For any small $\varepsilon>0$, if $a=a'=b=b'=\varepsilon$, then the tetrahedra between $\pp$ and $\mathcal{P'}$ admit positive dihedral angles satisfying the linear gluing equations (at all interior edges).
\end{prop}
\begin{proof}
Recall the Farey triangles $T_0, \dots,T_m$ separating $s$ from $s'$ (here, $m\geq 1$). By definition (see Section \ref{sec:fareytetrahedra}), we have $T_0=sAB$ and $T_m=s'A'B'$. Recall also the pleated surface $S_i$ associated to $T_i$: under our convention (transverse orientation for $S_i$), the pleating angles of $S_0$ at the edge pairs of slopes $A,B,s$ are $a,b,-a-b$ respectively. Similarly, the pleating angles of $S_m$ at the edge pairs of slopes $A',B',s'$ are $-a',-b',a'+b'$. We write these numbers in the corresponding corners of $T_0$ and $T_m$ (Figure \ref{fareyweights}).

For each $0<i<m$, the oriented line $\Lambda$ from $s$ to $s'$ enters $T_i$ across some Farey edge $e_i=T_i\cap T_{i-1}$, and exits through another edge $e_{i+1}$, either to the left or to the right: we say that $\Lambda$ \emph{makes a Left} or \emph{makes a Right} at $T_i$, and encode the combinatorics of $\Lambda$ into a word $\Omega=RLL...R$ of length $m-1$. 

No letter ($R$ or $L$) is associated \emph{a priori} to the Farey triangles $T_0$ and $T_m$. However, we will posit
that the path enters $T_0$ through the Farey edge $e_0:=sB$, and exits $T_m$ through the edge $e_{m+1}:=s'B'$, and associate the relevant letter ($R$ or $L$) to $T_0$ and to $T_m$. Hence, $\Omega$ becomes a word of length $m+1$. This convention is totally artificial, but it will allow us to streamline the notation in our argument.

\begin{figure}[h]
\psfrag{a}{$a$}
\psfrag{b}{$b$}
\psfrag{ab}{$\minus a\minus b$}
\psfrag{ap}{$\minus a'$}
\psfrag{bp}{$\minus b'$}
\psfrag{abp}{$a'\!\!\,+\!b'$}
\psfrag{aa}{$A$}
\psfrag{bb}{$B$}
\psfrag{app}{$A'$}
\psfrag{bbp}{$B'$}
\psfrag{s}{$s$}
\psfrag{sp}{$s'$}
\psfrag{leg}{($w_0=a~;~w_1=a+b$)}
\psfrag{w4}{$w_4$}
\psfrag{mw4}{$\minus w_4$}
\psfrag{w5}{$w_5$}
\psfrag{mw5}{$\minus w_5$}
\psfrag{ww}{$w_5\minus w_4$}
\psfrag{e0}{$e_0$}
\psfrag{e4}{$e_4$}
\psfrag{e5}{$e_5$}
\begin{center}
\includegraphics{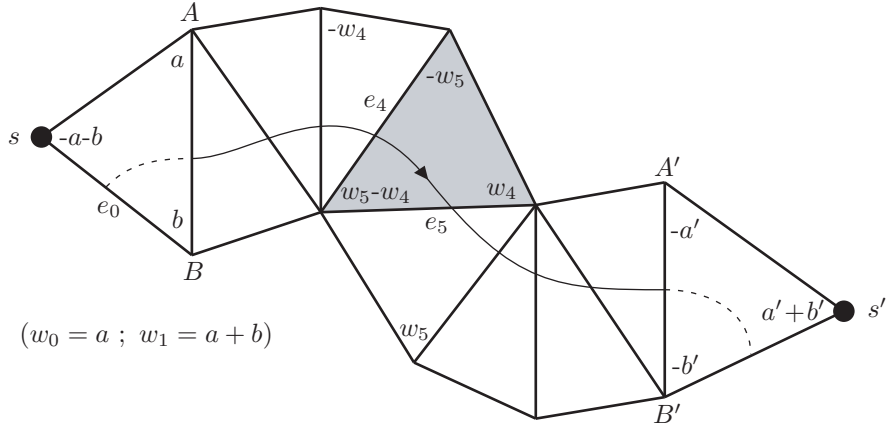}
\end{center}
\vspace{-1ex}
\caption{The pleating angles written in the corners of the Farey triangles $T_i$ associated to the pleated surfaces $S_i$.}
\label{fareyweights}
\end{figure}

For $1\leq i\leq m$, we denote by $\dot{\Delta}_i$ and $\ddot{\Delta}_i$ 
the two ideal tetrahedra separating the pleated surfaces $S_i$ and $S_{i-1}$: note that $\{\dot{\Delta}_i,\ddot{\Delta}_i\}$ is naturally associated to the Farey edge $e_i=T_i\cap T_{i-1}$. Our goal is to define dihedral angles for $\dot{\Delta}_i$ and $\ddot{\Delta}_i$ such that the linear gluing conditions around all edges are satisfied. We will in fact restrict to solutions invariant under the Klein group $V_4$, i.e. such that $\dot{\Delta}_i$ and $\ddot{\Delta}_i$ are isometric 
(\footnote{The graph carried by each $4$--punctured sphere $S_i$ is that of the edges of a tetrahedron, and its combinatorial symmetry group is $A_4$; the group $V_4 \subset A_4$ acts on these graphs in the usual way, by pairs of disjoint transpositions of ideal vertices.}) 
for all $i$ (this implies in particular that the angles of $\dot{\Delta}_i$ and $\ddot{\Delta}_i$ at any shared edge are equal). In what follows, $\Delta_i$ refers to any one of the ideal tetrahedra $\dot{\Delta}_i, \ddot{\Delta}_i$.

Denote by $\pi-w_i$ the dihedral angle of $\Delta_i$ at the pair of opposite edges that is not in $S_i\cap S_{i-1}$ (i.e. the pair of edges involved in the diagonal exchange that $\Delta_i$ represents). Then $S_i$ has one pleating angle equal to $w_i$ while $S_{i-1}$ has one pleating angle equal to $-w_i$. By translating indices, we find that for all $1\leq i \leq m-1$ the pleating angles of $S_i$ must be
$$-w_{i+1}~~~~~,~~~~~w_i~~~~~\text{ and }~~~~~ w_{i+1}-w_i$$
(the value of the third pleating angle is forced upon us by the condition that the pleating angles add up to $0$).
Further, we can write these three pleating angles in the corners of the Farey triangle $T_i$ associated to $S_i$ (this was partially done in Figure \ref{fareyweights}). In Figure \ref{fig:rl2} (top), 
denoting by $e_i$ the Farey edge $T_i\cap T_{i-1}$ associated to the tetrahedra $\{\dot{\Delta}_i,\ddot{\Delta}_i\}$, we see that $w_i$ is in the corner of $T_i$ opposite $e_i$, and $-w_{i+1}$ is in the corner of $T_i$ opposite $e_{i+1}$.

We repeat the same procedure for all indices $1\leq i\leq m-1$. It also extends naturally to $i=0$ and $i=m$ if we just set $(w_0,w_1)=(a,a+b)$ and $(w_m,w_{m+1})=(a'+b',a')$: we then recover the pleating angles of $S_0$ and $S_m$ defined previously.
 
The bottom part of Figure \ref{fig:rl2} shows the result 
of the labeling for two consecutive Farey triangles $T_{i-1}$ and  $T_i$, whose corresponding pleated surfaces $S_{i-1}$ and $S_i$ bound
the pair of isometric tetrahedra $\{\dot{\Delta}_i,\ddot{\Delta}_i\}$ (where $1\leq i \leq m$).
There are four possible cases, depending on the letters ($R$ or $L$) living on the Farey triangles $T_{i-1}$ and $T_i$. In order for $\Delta_i$ to have positive angles, assuming $0<w_i<\pi$ (for all $1\leq i \leq m$), it is necessary and sufficient that each pleating angle written just below the horizontal edge $e_i$ be larger than the pleating angle written just above, in Figure \ref{fig:rl2} (bottom): namely, the difference between these two pleating angles (of $S_i$ and $S_{i-1}$) is twice a dihedral angle of $\Delta_i$.

Suppose $(w_{i-1},w_i,w_{i+1})=(u,w,v)$.
Denoting by $x_i$ (resp. $y_i$) the angle of $\Delta_i$ at the edge whose slope is given by the right (resp. left) end of $e_i$ in Figure \ref{fig:rl2} (bottom), we thus find the following formulas for $x_i, y_i, z_i$:

\begin{equation} \label{interlettresarbo} \begin{array}{c|c|c|c|c}
\Omega &L~~~~~L         &R~~~~~R         &L~~~~~R        &R~~~~~L  \\ \hline
x_i    &\frac12(u+v)    &\frac12(-u+2w-v)&\frac12(u+w-v) &\frac12(-u+w+v)\\
y_i    &\frac12(-u+2w-v)&\frac12(u+v)    &\frac12(-u+w+v)&\frac12(u+w-v) \\
z_i    &\pi-w           &\pi-w           &\pi-w          &\pi-w \end{array}
\end{equation}

Define a \emph{hinge index} $i$ as an index such that the Farey triangles adjacent to the Farey edge $e_i$ carry different letters ($R$ and $L$). From (\ref{interlettresarbo}), we see that $\Delta_i$ has positive angles if and only if 
\begin{equation} \label{racohi} \left \{ \begin{array}{l}
\bullet \hspace{12pt} 0<w_i<\pi \text{ for all $1\leq i \leq m$ (range condition);}\\
\bullet \hspace{12pt} w_{i+1}+w_{i-1}<2w_i \text{ if $i$ is not a hinge index (concavity condition);}\\
\bullet \hspace{12pt} |w_{i+1}-w_{i-1}|<w_i \text{ if $i$ is a hinge index (hinge condition).} \end{array} \right . \end{equation}
Recall that $m\geq 1$. 
It is clear that there exist sequences $(w_0,\dots,w_{m+1})$ satisfying the above conditions such that $(w_0,w_1,w_m,w_{m+1})=(\varepsilon,2\varepsilon,2\varepsilon,\varepsilon)$: for example, set all $(w_i)_{1\leq i\leq m}$ equal to $2\varepsilon$, then perturb the non-hinge parameters among $\{w_i\}_{1<i<m}$ to obtain strong (piecewise) concavity. \end{proof}

\subsubsection{Trivial tangles} 

As in \cite[Appendix]{gf-bundle}, this method of constructing angle structures extends to the case when $\pp$ is glued to a trivial tangle ($1$--bracelet) of slope $s'$, realized by tetrahedra. Then, Figure \ref{fig:fold-clasp} on page \pageref{fig:fold-clasp} shows the surface $S_{m-1}$ associated to the next-to-last Farey triangle $T_{m-1}$. As in Figure \ref{fig:fold-clasp}, we assume (up to changing the marking of the $4$--punctured Conway sphere) that $s'=\frac{1}{2}$ and $T_{m-1}=1\infty 0$. Gluing the faces of $S_{m-1}$ in pairs amounts to requiring that $S_{m-1}$ have pleating angle $-\pi$ at the edges of slope $\infty$: in other words, that $w_m=\pi$. Thus, if we put $w_m=\pi$ (the Farey edge associated to $w_m$ is $e_m=T_m\cap T_{m-1}$), Table (\ref{interlettresarbo}) still allows us to derive the angles of the tetrahedra $\Delta_1,\dots,\Delta_{m-1}$, and the positivity conditions are still given by (\ref{racohi}). (In that case, there is no ``artificial'' letter living on $T_m$ and no ``artificial'' parameter $w_{m+1}$.) Angle structures are thus given by sequences \begin{equation} \label{eq:winailed}(w_0,w_1,w_2,\dots,w_{m-1},w_m)=(a,a+b,w_2,\dots,w_{m-1},\pi)\end{equation} which satisfy (\ref{racohi}) for $0<i<m$. To describe for which pairs $(a,b)$ such a sequence exists, we need some notation.

Given two rationals $q=y/x$ and $q'=y'/x'$ in reduced form
in $\mathbb{P}^1\mathbb{Q}$, define $$q\wedge q':=\left | \left | \begin{array}{cc}y&y'\\x&x'\end{array} \right | \right |$$ (absolute value of the determinant). We will use the following key property: if $pqr$ is a Farey triangle and $u,p,r,q$ are cyclically ordered in $\mathbb{P}^1\mathbb{Q}$, then $u\wedge r = (u\wedge p)+(u\wedge q)$. The property is easily checked in the case $(p,q)=(0,\infty)$ (where $u,r$ have opposite signs), and the general case follows because the $\wedge$-notation is invariant under $PSL_2(\mathbb{Z})$, which acts transitively on oriented Farey edges $pq$.

\begin{figure}[h!] \centering 
\psfrag{r}{$R$}
\psfrag{l}{$L$}
\psfrag{p}{$\pi$}
\psfrag{i}{$i$}
\psfrag{v}{$v_i$}
\psfrag{f}{$f(v_i)$}
\psfrag{h}{hinge}
\psfrag{o}{$0$}
\psfrag{1}{$1$}
\psfrag{m}{$m$}
\includegraphics{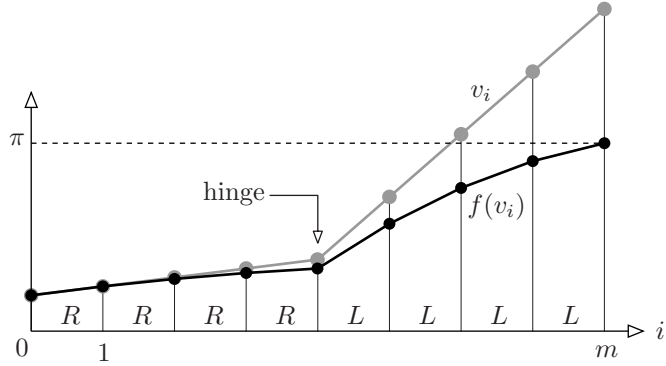}
\caption{The sequences $v=(v_i)_{0\leq i \leq m}$ and $w=f\circ v$.} \label{akashikaikyoo}
\end{figure}

\begin{prop} \label{prop:paqb}
Suppose a trivial tangle complement is glued to a large block $\pp$ that has non-negative pleating angles $a,b$ at the edge pairs of slope $A,B\in\mathbb{P}^1\mathbb{Q}$, satisfying (\ref{eq:aibi}). Suppose that $s,s' \in \mathbb{P}^1\mathbb{Q}$ are the preferred slopes of $\pp$ and of the trivial tangle, respectively. Then $sAB$ is the Farey triangle $T_0$; the points $s,A,s',B$ are cyclically ordered in $\mathbb{P}^1\mathbb{Q}$, and the tetrahedra $\Delta_1,\dots,\Delta_{m-1}$ (realizing the trivial tangle complement) admit positive angles if and only if
\begin{itemize}
\item $\hspace{20pt} s\wedge s'=2$ and $a+b=\pi$; or
\item $\hspace{20pt} s\wedge s'>2$ and $a(B\wedge s')+b(A\wedge s')>\pi>a+b$.
\end{itemize}
\end{prop}
\begin{proof}
The statements about the relative positions of $s,A,B,s'$ are true by construction and were proved in Section \ref{sec:fareytetrahedra}. The case $s\wedge s'=2$ corresponds to $m=1$, with the $4$--punctured sphere $S_{m-1}=S_0\subset \partial \pp$ being glued directly to itself (as in Figure \ref{fig:fold-clasp}). Since $m=1$, a sequence of the form (\ref{eq:winailed}) exists if and only if $a+b=\pi$. We now assume $s\wedge s'>2$, and consider the sequence of Farey triangles $T_0,\dots,T_m$ from $s$ to $s'$ (where $m\geq 2$). The inequality $\pi>a+b$ must clearly be true in (\ref{eq:winailed}) by (\ref{racohi}), so we focus on the other inequality (which says that $a,b$ are \emph{not too small}).

For each $0\leq i \leq m$, define $q_i$ to be the vertex of $T_i$ not belonging to the edge $e_i$ (where $e_0=sB$ and $e_i=T_i\cap T_{i-1}$ otherwise). In particular, $q_0=A$. If $$ \alpha_i=A\wedge q_i~~\text{ and }~~\beta_i=B\wedge q_i~~,$$
it is easy to check that both $(\alpha_i)$ and $(\beta_i)$ make the concavity and hinge conditions of (\ref{racohi}) \emph{critical} in the following sense: for each $0<i<m$,
\begin{itemize}
\item If $i$ is not a hinge index, then $\alpha_{i+1}+\alpha_{i-1}=2\alpha_i$ and $\beta_{i+1}+\beta_{i-1}=2\beta_i$;
\item If $i$ is a hinge index, then $\alpha_{i+1}=\alpha_i+\alpha_{i-1}$ and $\beta_{i+1}=\beta_i+\beta_{i-1}$.
\end{itemize}
(In the first case, observe that $\alpha_{i+1}-\alpha_i=A\wedge p=\alpha_i-\alpha_{i-1}$, where $p$ is the common vertex of the Farey edges $e_{i-1},e_i,e_{i+1}$. In the second case, observe that $q_{i-1}q_{i}q_{i+1}$ is a Farey triangle and $A,q_i,q_{i+1},q_{i-1}$ are cyclically ordered in $\mathbb{P}^1\mathbb{Q}$.) We say that $(\alpha_i)$ and $(\beta_i)$ satisfy \emph{the closure of} (\ref{racohi}) (the system obtained by turning all the strong inequalities of (\ref{racohi}) into weak ones).

Clearly, any linear combination of the sequences $(\alpha_i)$ and $(\beta_i)$ also makes the concavity and hinge conditions of (\ref{racohi}) critical. Define $$v_i:=a\beta_i+b\alpha_i=a(B\wedge q_i)+b(A\wedge q_i)~,$$ so that $(v_0,v_1)=(a,a+b)$, and $(v_i)$ satisfies the closure of (\ref{racohi}).

Note that $q_m=s'$, so $v_m=a(B\wedge s')+b(A\wedge s')$ is the left member of the inequality of the Proposition. 

\smallskip

\noindent {\bf Claim:} \emph{If $v'$ is another sequence which satisfies the closure of (\ref{racohi}) and $(v'_0,v'_1)=(v_0,v_1)$, then $v'_i\leq v_i$ and $v'_i-v'_{i-1} \leq v_i-v_{i-1}$ for all $1\leq i \leq m$.}

The claim is true for $i=1$, and follows in general by induction on $i$: if $i$ is not a hinge index, we have
$$\begin{array}{rcccccl}
v'_{i+1}&\leq&v'_i + (v'_i-v'_{i-1})&\leq&v_i + (v_i-v_{i-1})&=&v_{i+1}~; \\
v'_{i+1}-v'_i&\leq&v'_i - v'_{i-1}&\leq&v_i-v_{i-1}&=&v_{i+1} - v_i \end{array}$$
(in each line, the first inequality is true by (\ref{racohi}), and the second one by induction). Similarly, if $i$ is a hinge index, then
$$\begin{array}{rcccccl}
v'_{i+1}&\leq&v'_i + v'_{i-1}&\leq&v_i + v_{i-1}&=&v_{i+1}~; \\
v'_{i+1}-v'_i&\leq&v'_{i-1}&\leq& v_{i-1}&=&v_{i+1} - v_i~. \end{array}$$

Thus, if $v_m\leq\pi$, then no sequence $w=v'$ satisfies both (\ref{racohi}) and (\ref{eq:winailed}), so there can be no positive dihedral angle assignment for the tetrahedra $\Delta_1,\dots,\Delta_{m-1}$. Conversely, if $v_m>\pi$, we can define $w_i:=f(v_i)$ for all $1\leq i \leq n$, where \linebreak $f:[v_1,v_m]\rightarrow \mathbb{R}$ is increasing, strictly concave, $1$--Lipschitz and satisfies $f(v_1)=v_1=a+b$ and $f(v_m)=\pi$ (see Figure \ref{akashikaikyoo}); it is then straightforward to check that $(a=w_0,w_1,\dots,w_m=\pi)$ satisfies (\ref{racohi}) --- except of course the range condition at $w_m=\pi$. \end{proof}

\subsubsection{2-bridge links}

When two trivial tangles are glued together, we obtain a $2$-bridge link. The construction preceding Proposition \ref{prop:paqb} can be performed both near $s$ and near $s'$, and we refer to \cite{gf-bundle} for a much more complete treatment. In the remainder of the paper, we assume that the candidate link $K$ is not a $2$-bridge link, i.e. contains at least one large bracelet.

\subsection{Non (strongly) Montesinos links}\label{sec:non-montesinos}

Observe that the condition in each of the Propositions \ref{prop:chi-angles}, \ref{prop:degreed}, \ref{prop:degree3} requires that the angles $a_i, b_i$ be, in a loose sense, \emph{small enough}, while Proposition \ref{prop:paqb} requires them to be \emph{large enough}. The conflict that arises causes some arborescent link complements to have no angle structures (essentially, the third family of exceptions in Theorem \ref{thm:main}). 

\begin{define}
A \emph{strongly Montesinos} link is an arborescent link which, after the reduction of Section \ref{sec:algorithm}, consists of one non-augmented $d$--bracelet ($d\geq 3$) attached to $d$ trivial tangles.
\end{define}

Note that some very simple Montesinos links are not strongly Montesinos: for example, links with tangle slopes $(\pm \frac{1}{2}, \pm \frac{1}{2}, \frac{y}{x})$ were reduced in step $(5)$ of the algorithm of Section \ref{sec:algorithm}. Such links consist of an augmented $1$--bracelet glued to a non-augmented $1$--bracelet, and may or may not be candidate links, depending on whether the preferred slopes of the two bracelets satisfy the minimum--distance table of Proposition \ref{prop:algorithm-output}.

\begin{prop}
If the candidate link $K$ is not strongly Montesinos, then we can endow all blocks with non-negative dihedral angles satisfying (\ref{eq:aibi}) such that the girth condition (\ref{eq:chi-angles}) holds for all non-augmented blocks, and the condition of Proposition \ref{prop:paqb} holds at each trivial tangle. As a result, the ideal decomposition of the link complement admits angle structures, and the link is hyperbolic.
\end{prop}
\begin{proof}
Consider a non-augmented block $\pp$: since the link is not Montesinos, 
$\pp$ is separated by a product region from some other large block (augmented or not). By Propositions \ref{prop:empty-branch}--\ref{prop:innerbranch}, the dihedral angles $a_k, b_k$ of $\pp$ at the boundary of that product region can be taken smaller than or equal to any small $\varepsilon_0>0$, and the product region will still admit positive angle structures. More precisely, we take $a_k=b_k=\varepsilon_0$ (Prop. \ref{prop:innerbranch}) except in the special case where the two blocks are glued directly to one another: then, we use the observation that immediately follows Proposition \ref{prop:empty-branch} and take $a_k, b_k$ equal to $\varepsilon_0$ and $0$, though not necessarily in that order.)

We must now find dihedral angles for $\pp$ such that the girth condition (\ref{eq:chi-angles}) holds. Consider a trivial tangle attached to $\pp$, along the band $U_i$. 
Note that 
$$A\wedge s' + B\wedge s' \; = \; s\wedge s' \; \geq \; 2,$$ 
where the inequality follows from the table of minimal Farey distances in Proposition \ref{prop:algorithm-output}.
Therefore, by Proposition \ref{prop:paqb}, the tetrahedra in the trivial tangle will admit positive structures as soon as 
\begin{equation}\label{eq:balance}
a_i=b_i=\left\{ 
\begin{array}{cl}
\displaystyle{\frac{\pi}{2}}~ & \mbox{if } s \wedge s' = 2, \\
\varepsilon +\displaystyle{\frac{\pi}{s\wedge s'}} & \mbox{if } s \wedge s' > 2,
\end{array}
\right.
\end{equation}
for a small $\varepsilon>0$.
Since $\text{max}\,\{a_k, b_k\}=\varepsilon_0$, under this choice of values, the right member $\sum_{i=1}^d \max\{a_i,b_i\}$ of the girth condition (\ref{eq:chi-angles}) is thus at most $$\varepsilon_0+(d-1)\frac{\pi}{2}~.$$ If $d\geq 4$, this quantity is already less than the left member $(d-2)\pi$ of (\ref{eq:chi-angles}). If $d=3$, recall from Section \ref{sec:algorithm} that $\pp$ is not attached to two trivial tangles of slope $1/2$ (otherwise, we would have replaced the $3$--bracelet associated to $\pp$ by an augmented $1$--bracelet in step $(5)$ of the algorithm). Thus, the above upper bound can be further improved to $\varepsilon_0+(\pi/2)+(\varepsilon+\pi/3)<(d-2)\pi=\pi$, so the girth condition (\ref{eq:chi-angles}) is satisfied. 

As for augmented bracelets, there is nothing to check: as soon as the $a_i, b_i$ for an augmented block $\mathcal{P}$ satisfy (\ref{eq:aibi}) and (near trivial tangles) the condition of Proposition \ref{prop:paqb}, $\mathcal{P}$ is realized by a union of tetrahedra with positive angles, by Equation (\ref{eq:aug-angles}) above.

Finally, Theorem \ref{thm:block-hyperbolic} implies the existence of a hyperbolic structure.
\end{proof}

\subsection{Strongly Montesinos links}\label{sec:montesinos}

Suppose the candidate link $K$ is strongly Montesinos, and recall the non-negative integers $n_{\mathcal{A}},n_{\mathcal{B}}$ from 
Definition \ref{def:na-nb}.

\begin{prop}
If $n_{\mathcal{A}}\geq 3$ or $n_{\mathcal{B}}\geq 3$, the block decomposition admits angle structures.
\end{prop}
\begin{proof}
Assume $n_{\mathcal{A}}\geq 3$ (it is enough to treat this case). It is straightforward to find a pair $(\alpha, \beta)$ satisfying the condition of Proposition \ref{prop:degreed} (namely $\alpha \geq \beta$), such that $(a,b)=(\pi-\alpha\varepsilon, \beta \varepsilon)$ satisfies the condition of Proposition \ref{prop:paqb} for small $\varepsilon$. For example, take $\beta=\alpha$ if $s\wedge s'=2$, and $\beta=\frac{2}{3}\alpha$ if $s\wedge s'\geq 3$ (recall $(A\wedge s')+(B\wedge s')=s\wedge s'$,
so $a(B\wedge s')+b(A\wedge s')\geq \text{min}\,\{2a+b,a+2b\}>\pi$). The conclusion now follows from the two quoted propositions.
\end{proof}

The next two Propositions deal exactly with the remaining strongly Montesinos links, where $\max \{ n_{\mathcal{A}}, n_{\mathcal{B}} \}=2$. In each Proposition, we find a few non-hyperbolic links: 
\begin{itemize}
\item In Proposition \ref{prop:chain-link}, it is the link already mentioned in Figure \ref{fig:chain-link} and Remark \ref{rem:chain-link}, which falls into the second class of exceptions of Theorem \ref{thm:main}.
\item In Proposition \ref{prop:realwork}, it is exactly the strongly Montesinos links among the third class of exceptions of Theorem \ref{thm:main}.
\end{itemize}

Recall (Definition \ref{def:na-nb} and Remark \ref{rem:count-twists}) that $n_{\mathcal{B}}$ is the number of twists in the ``main band'' of 
a strongly Montesinos link $K$,
as defined prior to Definition \ref{def:candidate}: thus, it would be a straightforward exercise to translate the current block presentations back into planar link diagrams. 
Recall as well that we have chosen markings for the Conway spheres of $K$, in which the large bracelet has preferred slope $\infty$ and the trivial tangles have preferred slopes in the interval $(0,1)$.
All the exceptions arising in Propositions \ref{prop:chain-link}--\ref{prop:realwork} were preemptively ruled out by the last condition in the definition \ref{def:candidate} of candidate links.

\begin{prop} \label{prop:chain-link}
If $d=4$ and $n_{\mathcal{A}}=n_{\mathcal{B}}=2$, the block decomposition admits angle structures, unless all the trivial tangles have slope $1/2$.
\end{prop}
\begin{proof}
The quadruple $(s,A,B,s')$ associated to a trivial tangle of slope $\frac{y}{x}$ is by definition $(\infty,1,0,\frac{y}{x})$: therefore, the key condition $a(s'\wedge B)+b(s'\wedge A)>\pi$ from Proposition \ref{prop:paqb} becomes \begin{equation} \label{eq:paqbtangle} a(y)+b(x-y)>\pi~.\end{equation} Moreover, the denominator $x$ of the slope $\frac{y}{x}$ of the trivial tangle is the integer $\infty \wedge \frac{y}{x}=s\wedge s'$. If at least one of these denominators is larger than $2$, we can set $a_i, b_i$ as in (\ref{eq:balance}) above, and immediately obtain the girth condition (\ref{eq:chi-angles}) because $\frac{\pi}{2}+\frac{\pi}{2}+\frac{\pi}{2}+(\varepsilon+\frac{\pi}{3})<2\pi$. If all denominators are $2$, the link is not hyperbolic (Figure \ref{fig:chain-link}) and not candidate; it belongs to the second family of exceptions of Theorem \ref{thm:main}.
\end{proof}

\begin{prop} \label{prop:realwork}
If $d=3$ and $(n_{\mathcal{A}}, n_{\mathcal{B}})=(2,1)$, assume the trivial tangles have slopes $\frac{y_1}{x_1},\frac{y_2}{x_2},\frac{y_3}{x_3} \in (0,1)$: the block decomposition admits angle structures, unless one has $y_1=y_2=y_3=1$ and $\frac{1}{x_1}+\frac{1}{x_2}+\frac{1}{x_3}\geq 1$.
\end{prop}
\begin{proof}
Again, if $\sum \frac{1}{x_i}<1$, we can set $a_i, b_i$ as in (\ref{eq:balance}) to obtain the girth condition (\ref{eq:chi-angles}). Thus, assume $\sum \frac{1}{x_i}\geq 1$ and (up to a permutation) $y_3\geq 2$. This entails in particular that $x_3\geq 3$.

We will set $(a_i, b_i)=(\pi-\alpha_i, \beta_i)$ for well-chosen \emph{small} positive numbers $\alpha_i, \beta_i$. As in (\ref{eq:paqbtangle}) above, the key condition from Proposition \ref{prop:paqb} is still $a_i y_i+b_i(x_i-y_i)>\pi$. If $y_i>1$, this condition is vacuous for small $\alpha_i,\beta_i$. If $y_i=1$, it can be written $(x_i-1)\beta_i>\alpha_i$. Thus, the full set of 
sufficient conditions to be satisfied is:

\begin{itemize}
\item If $x_i=2$ then $\alpha_i=\beta_i>0$ (see Proposition \ref{prop:paqb});
\item If $x_i>2$ then $\alpha_i>\beta_i>0$ (see Proposition \ref{prop:paqb});
\item If $x_i>2$ and $y_i=1$ then $(x_i-1)\beta_i>\alpha_i$ (see Proposition \ref{prop:paqb});
\item For all $i\in \{1,2,3\}$, we have $\alpha_i>\beta_{i+1}+\beta_{i-1}$, taking indices modulo $3$ (see Proposition \ref{prop:degree3}).
\end{itemize}
(The first three conditions ensure the existence of angle structures for the trivial tangles; the last one, for the solid torus).

If $x_1=x_2=x_3=3$, we take $$\begin{array}{rrrclrrrc}
(~\beta_1~,&\beta_2~,&\beta_3~)&
=&(&\varepsilon~,&\varepsilon~,&\mu\varepsilon~)&\\
(~\alpha_1~,&\alpha_2~,&\alpha_3~)&=&(&(2-\mu)\varepsilon~,&(2-\mu)\varepsilon~,&M\varepsilon~)&;\end{array}$$ the conditions above are clearly satisfied if the positive parameters $\mu,M$ verify $\mu<\frac12$ and $M>2$. For example, $(\mu,M)=(\frac{1}{3},3)$.

Finally, if $x_1=2$, then $x_2\geq 3$: otherwise, the $3$--bracelet associated to the block $\pp$ would have been replaced by an augmented $1$--bracelet in step $(5)$ of the algorithm of Section \ref{sec:algorithm} (so the link would not be strongly Montesinos). We thus set
$$\begin{array}{rrrclrrrc}
(~\beta_1~,&\beta_2~,&\beta_3~)&=&(&(1+m)\varepsilon~,&\varepsilon~,& \mu \varepsilon~)&\\
(~\alpha_1~,&\alpha_2~,&\alpha_3~)&=&(&(1+m)\varepsilon~,&(x_2-1-\mu)\varepsilon~,&M\varepsilon~)&;\end{array}$$
the conditions above are clearly satisfied if the positive parameters $\mu,m,M$ satisfy $\mu<m$ and $2\mu+m<1$ and $M>2+m$. For instance, $(\mu,m,M)=(\frac{1}{4},\frac{1}{3},3).$ 
\end{proof}

Thus, all candidate links are hyperbolic. We have proved Theorem \ref{thm:sub-main}, hence Theorem \ref{thm:main}.

\bibliographystyle{hamsplain}

\bibliography{biblio.bib}
\end{document}